\documentclass[10pt,final,twocolumn,twoside]{IEEEtran}

\newcommand{\DRAFT}[1]{}

\usepackage{amssymb,amsmath}

\usepackage{graphicx}
\usepackage{algorithmic}
\usepackage{setspace}
\graphicspath{{images/}}

\newenvironment{changemargin}[2]{\begin{list}{}{%
\setlength{\topsep}{0pt}%
\setlength{\leftmargin}{0pt}%
\setlength{\rightmargin}{0pt}%
\setlength{\listparindent}{\parindent}%
\setlength{\itemindent}{\parindent}%
\setlength{\parsep}{0pt plus 1pt}%
\addtolength{\leftmargin}{#1}%
\addtolength{\rightmargin}{#2}%
}\item }{\end{list}}

\usepackage[latin1]{inputenc}
\usepackage{citesort}
\usepackage{xspace}
\usepackage{url}

\input{alphabet.tex}

\newsavebox{\fminibox}
\newlength{\fminilength}

\newtheorem{remark}{Remark}

\newcommand{\diag}[1]{\ensuremath{\mathrm{diag}\{#1\}}\xspace}

\def\wrt{with respect to\xspace}
\def\cc#1{\setlength{\tabcolsep}{0pt}\begin{tabular}{c}#1\end{tabular}}

\newcommand{\Go}{\ensuremath{\Gv_\oD}\xspace}
\def\T{^\tD}
\renewcommand{\Gc}{\ensuremath{\Gv_\cD}\xspace}

\newcommand{\x}{\ensuremath{\xb}\xspace}

\def\Px{\Pv_{\kern-2pt\xb}}
\def\Pw{\Pv_{\kern-2pt\wb}}
\def\argmin{\mathop{\mathrm{argmin}}} 

\def\pth#1{\left(#1\right)}                
\def\acc#1{\left\{#1\right\}}

\def\norm#1{\left\|#1\right\|}             \def\stdnorm#1{\|#1\|}

\def\ldotsv{,\,\ldots,\,}
\def\reg{_{\mathrm{r}}}
\def\CSI{_{\mathrm{CSI}}}
\def\ie{\textit{i.e.},\XS}
\def\rem#1{}                    
\def\froc#1#2{{#1/#2}}
\def\old{^\text{old}}

\usepackage{color,soul}

\definecolor{lightblue}{rgb}{.4,.9,1}

\definecolor{orange}{rgb}{1,.9,.5}

\begin{document}

\title{On Algorithms Based on Joint Estimation of Currents and Contrast in Microwave Tomography}

\author{%
Paul-André Barrière, Jérôme Idier, Yves Goussard, and Jean-Jacques Laurin
\thanks{Paul-André Barrière is with \'Ecole Polytechnique de Montréal, Dep. of
Electrical Engineering, CP 6079, Succ. Centre-ville Montréal (Qu\'ebec), H3C~3A7, Canada
and with Institut de Recherche en
Communications et en Cybernétique de Nantes (IRCCyN), 1 rue de la Noë, BP 92101, 44321 Nantes Cedex 03
(email: paul-andre.barriere@polymtl.ca).}
\thanks{Jérôme Idier is with Institut de Recherche en
Communications et en Cybernétique de Nantes (IRCCyN), 1 rue de la Noë, BP 92101, 44321 Nantes Cedex 03 (email:
Jerome.Idier@irccyn.ec-nantes.fr).}
\thanks{Yves Goussard and Jean-Jacques Laurin are with \'Ecole Polytechnique de Montréal, Dep. of Electrical
Engineering, CP 6079, Succ. Centre-ville Montréal (Qu\'ebec), H3C~3A7, Canada
(email: yves.goussard@polymtl.ca; jean-jacques.laurin@polymtl.ca).}}

\maketitle
\begin{abstract}
This paper deals with improvements to the contrast source inversion
method which is widely used in microwave tomography. First, the
method is reviewed and weaknesses of both the criterion form and the
optimization strategy are underlined. Then, two new algorithms are
proposed. Both of them are based on the same criterion, similar but
more robust than the one used in contrast source inversion. The
first technique keeps the main characteristics of the contrast
source inversion optimization scheme but is based on a better
exploitation of the conjugate gradient algorithm. The second
technique is based on a preconditioned conjugate gradient algorithm
and performs simultaneous updates of sets of unknowns that are normally
processed sequentially. Both techniques are shown to be more efficient
than original contrast source inversion.
\end{abstract}


\section{Introduction}
\label{sec:intro} The objective of microwave tomography is to reconstruct the permittivity and conductivity distributions of an object
under test. This is performed from measurements of the field scattered by this object under various conditions of illumination. Microwave tomography has shown
great potential in several application areas, notably biomedical imaging, non destructive testing and geoscience.

Unlike other well-known imaging techniques (e.g. X-ray tomography), the involved wavelengths are long compared to the structural features
of the object under test. Consequently, ray propagation approximation is not suitable. Instead an integral equation formulation, which is nonlinear and ill-posed \cite{Isernia_97_NE}, must be used.

A large number of techniques have been proposed to perform inversion (\ie finding the permittivity and conductivity distributions from the
measurements). In most of them inversion is formulated as an optimization problem. Differences between methods then come from the nature of the
criterion and the type of minimization algorithm.

Proposed criteria are made up of one, two or three terms with the possible addition of a regularization term. When only one term is used, the
quadratic error between actual and predicted measurements is minimized. In the two- and three-term cases, criteria penalize both the error on
the measurements and the error on some constraint inside the domain of interest.

The problem being nonlinear, the criteria to be minimized may have local minima. To avoid being trapped in them, some methods use global
minimization algorithms~\cite{Benedetti_05_EE,Donelli_05_CA,Massa_05_PG,Caorsi_04_MM,Donelli_05_CA,Garnero_91_MI}. However, as the complexity of
such algorithms grows very rapidly with the number of unknowns, most proposed methods use local optimization schemes and assume that no local
minima will be encountered.

Some of these methods are based on successive linearizations of the problem using the Born approximation or some
variants~\cite{Semenov_99_TD,van_den_Berg_02_IS,Liu_02_AM,Abubakar_04_RI,Souvorov_98_MT,Franchois_97_MI,Joachimowicz_98_CS,Joachimowicz_91_IS,Bandyopadhyay_05_IA}.
The capacity of these methods to provide accurate solutions to problems with large scatterers varies according to the used approximation. For instance, the distorted Born method~\cite{Chew_90_RT} (which is equivalent to the Newton-Kantorovich~\cite{Joachimowicz_91_IS,Remis_00_EN} and
Levenberg-Marquardt~\cite{Franchois_97_MI} methods) is more efficient for large scatterers than the Born iterative method~\cite{Li_04_TD} or
some extensions like the one presented in~\cite{Song_04_FT}. Nevertheless, most of these methods require complete computation of the forward
problem at each iteration which is computationally burdensome.

To avoid the necessity for solving the forward problem at each iteration, two-term criterion techniques were proposed. The idea is to performed
optimization not only with respect to electrical properties but also \wrt the total field in the region of interest. This is done with a
so-called modified gradient technique in \cite{Kleinman_93_ER}, a conjugate gradient algorithm in \cite{Barkeshli_94_IM} and a quasi-Newton
minimization procedure in \cite{Isernia_97_NE}.

The contrast source inversion (CSI) method \cite{van_den_Berg_97_CS,van_den_Berg_02_IS,Abubakar_01_TV,Abubakar_02_IB,Omrane_06_SR} is also a 2-term criterion technique but this time the problem is formulated as a function of equivalent currents and electrical properties of the object under test (instead of the electric field and electrical properties). Optimization is performed by a conjugate gradient algorithm with an alternate update of the unknowns.

We mentioned that, due to the nonlinearity of the problem, global optimization schemes could be a good choice. In
\cite{Benedetti_05_EE,Donelli_05_CA,Massa_05_PG,Caorsi_04_MM,Donelli_05_CA}, such algorithms are proposed for optimization of a two-term
criterion. In \cite{Garnero_91_MI} a simulated annealing algorithm is used to minimize a one-term criterion.

In this paper, the emphasis is placed on well-known CSI methods which offer a good compromise between quality of the solution and
computational effort. We address both questions of the form of the criterion and of the corresponding optimization algorithms. We first give
background information on the CSI method; then we show that the corresponding criterion exhibits unwanted characteristics and could lead to
degenerate solutions. Two pitfalls of the optimization scheme are also underlined.

Based on this analysis two new methods are proposed; both use the same criterion, the form of which is deduced from optimization and
regularization theories. Therefore, these techniques only differ by their respective minimization strategies. The first one retains the main
characteristics of the CSI method but is based on a better exploitation of the linear conjugate gradient algorithm. The second one is based on
simultaneous updates of sets of unknowns and a preconditioned nonlinear conjugate gradient algorithm is used to perform the minimization. The
two proposed methods exhibit improved robustness and convergence speed compared to CSI.

Note that the scope of the study is limited to analysis of and improvements to techniques that retain the main characteristics of the original CSI methods. Further performance gains, either by significantly altering the nature of the objective function or by using approximate forms thereof, are investigated in~\cite{Barriere_08_CF} and~\cite{Barriere_08_NA}.


\section{Context}
\label{sec:context} 
\begin{figure}[htpb]
\centerline{\includegraphics[width=4cm]{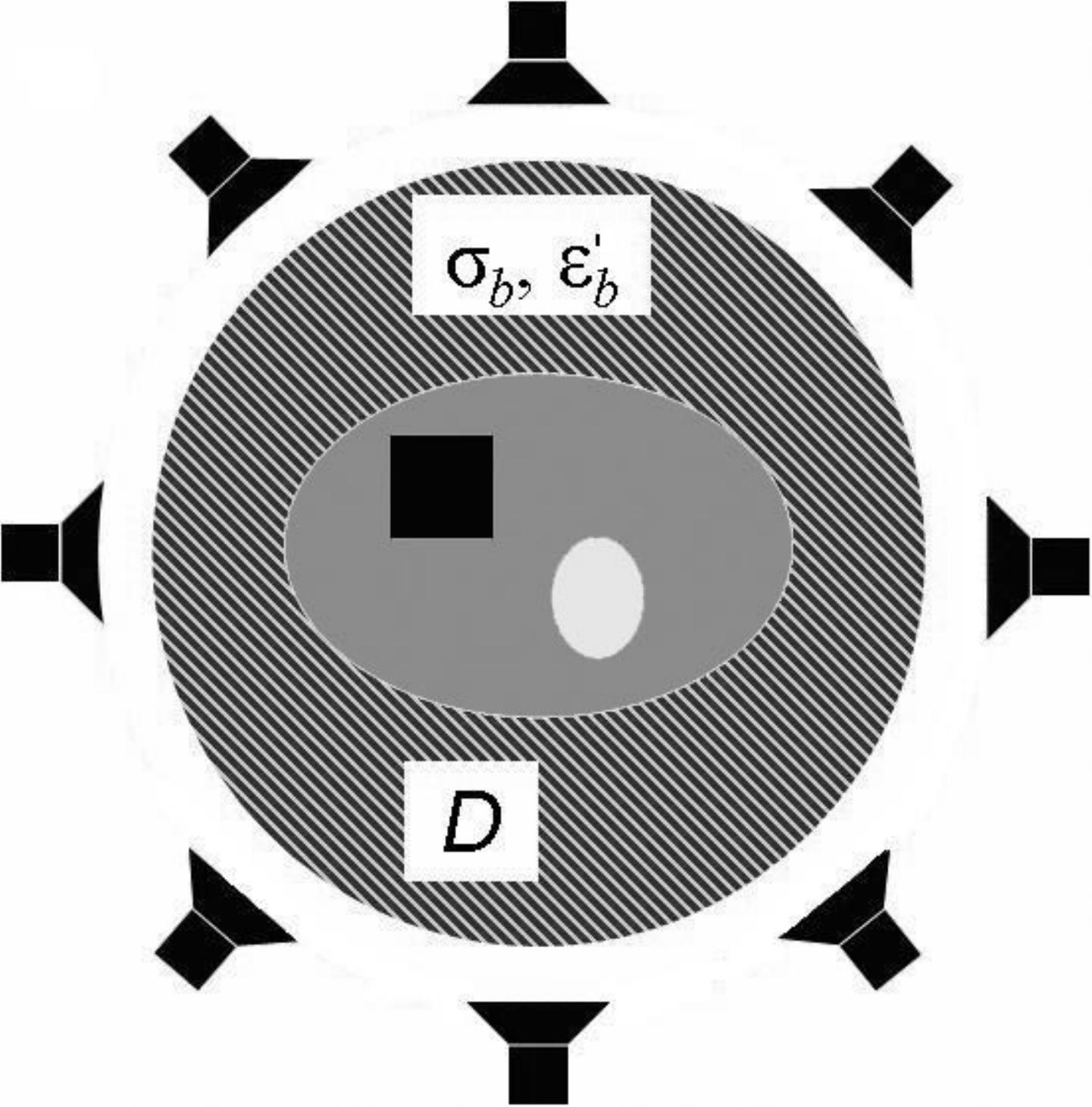}}
\caption{Typical microwave tomography setup} \label{fig:set_up}
\end{figure}

In a microwave tomography experiment the objective is to find permittivity ($\epsilon'$) and conductivity ($\sigma$) distributions of an object under test
placed into a domain $D$ and sequentially illuminated by $M$ different microwave emitters. For each illumination, the scattered field is
gathered at $N$ points. A typical setup is depicted in Fig.~\ref{fig:set_up}. In this study, we limit ourselves to the 2-D TM case. This means
we assume that all quantities are constant along the $z$ direction and that the electric field is parallel to it.

By refering to the volume equivalence theorem and by using the method of moments as a discretization technique, the microwave tomography experiment can be
described by the following system of equations (for more detail on the derivation we refer to \cite{Abubakar_02_IB}):

For $i=1\ldotsv M$,
\begin{subequations}
\label{eq:TMO_CS}
\begin{equation}
\label{eq:TMO_CS_obs} \yb_i=\Gv_\oD\wb_i + \nb_i^\oD
\end{equation}
\begin{equation}
\label{eq:TMO_CS_champs} \wb_i=\Xv(\Eb^0_i+\Gv_\cD\wb_i)
\end{equation}
\end{subequations}
where bold-italic fonts and bold-straight fonts denote column vectors and matrices, respectively. $\yb_i$ is a length $N$ vector that contains
the measured scattered field related to the $i$th illumination. $\Eb^0_i$ is the discretized incident field in $D$ (\ie the field when $D$ is
completely filled with the background medium). Its length is the number of discretization points, denoted by $n$. $\Gv_\oD$ and $\Gv_\cD$ are
Green matrices of size $N \times n$ and $n \times n$, respectively. $\Xv$ is a diagonal matrix such that $\Xv=\diag{\xb}$. $\xb$ and $\wb_i$ are
vectors of length $n$ and are called contrast vector and current vector, respectivley. Finally, $\nb_i^\oD$ is a noise vector that models all
perturbations encountered in a microwave tomography experiment.

Equations~\eqref{eq:TMO_CS_obs} and~\eqref{eq:TMO_CS_champs} are called \textit{observation} and \textit{coupling} equations, respectively. The
unknown quantities are $\xb$ and $\Wv=(\wb_1\ldotsv\wb_M)$ while \xb
represents the actual quantity of interest, since it contains all relevant
information about the permittivity and conductivity of the object under test. Indeed, $\xb$ is related to the electrical characteristics as follow:
\begin{equation}
\label{eq:def_x} \xb=(\boldsymbol\epsilon^d-\boldsymbol\epsilon^d_b)/\boldsymbol\epsilon^d_b
\end{equation}
where the division is performed term by term and where $\boldsymbol\epsilon^d$ and $\boldsymbol\epsilon^d_b$ are the discretized complex permittivity of the object under test and of background medium, respectively. Complex
permittivity is expressed by
\begin{equation}
\label{eq:perm_complexe} \epsilon=\epsilon'-j\sigma/\omega
\end{equation}
where $\omega$ denotes the angular frequency. Estimating $\xb$ is then equivalent to estimate $\epsilon$ and $\sigma$.

It is possible to eliminate $\wb_i$ from~\eqref{eq:TMO_CS}. The following equation is then obtained:
\begin{equation}
\label{eq:TMO_non_lin} \yb_i=\Gv_\oD(\Iv-\Xv\Gv_\cD)^{-1}\Xv\Eb^0_i +\nb_i^\oD
\end{equation}
where $\Iv$ is the identity matrix. Note that now, $\yb_i$ is a function of $\xb$ only. The above expression highlights the nonlinearity of the
problem which greatly complicates its resolution.

Some approaches are based on (or equivalent to) solving~\eqref{eq:TMO_non_lin} in the least-squares
sense~\cite{Carfantan97,Franchois_97_MI}. Meanwhile, those methods necessitate a high calculation cost. In order to circumvent the difficulty, two term criterion methods are based on system \eqref{eq:TMO_CS} (or an equivalent system
based on the total field and the contrast instead of the currents and the contrast \cite{Abubakar_02_IB,Kleinman_93_ER}) and minimize the sum of
errors on both observation and coupling equations. Such an approach is not equivalent to solving \eqref{eq:TMO_non_lin}, since the coupling
equation is not fulfilled exactly. However, under most circumstances, the solutions are extremely similar. Therefore, the rest of the study is
dedicated to this type of approach.


\section{Analysis of the CSI method }
In this section we present the background results on CSI and underline some of weaknesses. Both the criterion form and the optimization scheme
are studied.

\subsection{Background results on CSI}
\label{sec:analyse_crtiere} CSI formulates the microwave tomography problem in the framework of optimization. The objective is to estimate
$(\hat{\xb},\hat{\Wv})$ minimizing  a given criterion  $F$. This criterion makes a trade-off between the errors that respectively
affect the observation \eqref{eq:TMO_CS_obs} and coupling \eqref{eq:TMO_CS_champs} equations. A possible regularization term can be added to deal with the ill-posed nature of the problem~\cite{Isernia_97_NE}. More precisely, we have
\begin{equation}
\label{eq:F_sum_de_F} (\hat{\xb},\hat{\Wv})=\argmin_{\xb,\Wv}F
\end{equation}
\begin{subequations}
\begin{align}
\label{eq:2termes_typique}
&F=F_1+\lambda F_2+\lambda\reg F\reg\\
\label{eq:F1}
&F_1=\sum_i\|\yb_i-\Gv_\oD\wb_i\|^2 \\
\label{eq:F2}
&F_2=\sum_i\|\Xv(\Eb^0_i+\Gv_\cD\wb_i)-\wb_i\|^2\\
\label{eq:Freg} &F\reg=\phi(\xb)
\end{align}
\end{subequations}
where $\lambda$ and $\lambda\reg$ respectively denote the weight and the regularization factor, $\norm{\cdot}$ is the Euclidean norm, and $\phi$
is a regularization function.

\begin{changemargin}{0.33cm}{0cm}
\remark{In~\cite{Abubakar_01_TV}, the concept of multiplicative regularization was introduced. The idea is to multiply $F_1 + \lambda F_2$ by
the penalty term $F\reg$, instead of adding it. Here, we shall not pursue in this direction; we shall rather adopt the classical penalized
least-square framework.}
\end{changemargin}

CSI proposes the following explicit expression of $\lambda$ as a function of \xb:
\begin{equation}
\label{eq:lambda_CSI}
\lambda=\lambda\CSI=\frac{\sum_i\|\yb_i\|^2}{\sum_i\|\Xv\Eb^0_i\|^2}\ \ .
\end{equation}
This choice is justified by the fact that $F_1$ and $\lambda F_2$ are equal when the currents vanish.

The optimization is performed according to a block-component scheme based on alternate updates of each set of unknowns. More precisely,
minimization \wrt each $\wb_i$ is performed for fixed values of $\xb$ and $\wb_j$, $j\neq i$; then, $\Wv$ is fixed at its current value and
minimization is performed \wrt $\xb$. These updates form one iteration of the algorithm, and they are repeated until convergence. Each update of
the CSI method consists of a single conjugate gradient step. Table~\ref{algo:CSI} presents the algorithm. More details are given in
\cite{Abubakar_01_TV}.

\begin{table}[t]
\setlength\fboxsep{2mm}
\begin{spacing}{1.2}
\fbox{\kern-8pt\begin{minipage}{.985\linewidth}
\smallskip
\begin{algorithmic}
\STATE $k \gets 0$ \STATE Initialize \xb and $\Wv$ \REPEAT \FOR{$i=1,...,M$} \STATE Perform one iteration of the conjugate gradient algorithm to minimize $F$
\wrt $\wb_i$ \ENDFOR \STATE Perform one iteration of the conjugate gradient algorithm to minimize $F$ \wrt $\xb$ \STATE $k \gets k+1$ \UNTIL{Convergence}
\end{algorithmic}
\end{minipage}}
\end{spacing}
\caption{Implementation of the CSI method according to \cite{Abubakar_01_TV}} \label{algo:CSI}
\end{table}

\subsection{Pitfalls of the CSI method}
\label{subsec:weak} It seems that both the criterion and the optimization algorithm proposed in the CSI method suffer from some weaknesses.

\subsubsection{Criterion}
\label{subsubsec:criterion} The choice \eqref{eq:lambda_CSI} for $\lambda$ seems to be inappropriate essentially because the resulting criterion
$F$ reaches its minimum value for degenerate solutions if no regularization is used. More precisely, it is shown in Appendix~\ref{degenerate}
that solutions ($\hat{\xb}$,$\hat{\Wv}$) exist such that $\|\hat{\xb}\|\rightarrow\infty$ and $F_1 + \lambda\CSI F_2\rightarrow0$.  This is
obviously an undesirable feature, since it shows that any globally converging minimization method will produce a degenerate solution. As for
more realistic, locally converging methods, they could either converge towards a local minimum, possibly near the expected solution, or towards
a degenerate global minimum.  In Section~\ref{results}, an example is shown where a local optimization algorithm converges toward a
degenerate solution.

From a theoretical standpoint, an additional regularization term could solve the problem, provided that $F\reg(\xb)\to\infty$ when
$\norm{\xb}\to\infty$.  However, local (or even global) minima with large norms may still exist if $\lambda\reg$ is not chosen large enough,
while too large values of $\lambda\reg$ will produce overregularized solutions.

From a practical standpoint, CSI is widely used and, to our knowledge, convergence to degenerate solutions has not been reported yet. A reason
could be that, in a widespread version of the method, the dependence of $\lambda\CSI$ on \xb is omitted in the calculation of the gradient
component $\nabla_\xb F$.
However, it is shown in the next section that this leads to undesirable properties of the optimization algorithm.

Finally, it should be underlined that the choice \eqref{eq:lambda_CSI} is not justified by solution quality or computation cost arguments. It
could then be expected that better choices be possible from those points of view.

\subsubsection{Optimization scheme}
\label{subsubsec:CSI_algo_weak}

The first problem with the CSI optimization scheme is to rely on conjugacy formulas in an unfounded way. Indeed, the conjugate gradient algorithm has been
designed to minimize a criterion \wrt one set of unknowns \cite{Nocedal_livre_NO}. The conjugacy of the descent directions can be defined only
if this criterion remains the same during the optimization process. Meanwhile, in CSI, the criterion \wrt each set of unknowns change between
two updates. For instance, between two updates of $\wb_i$ at iterations $k$ and $k+1$, the criterion $F(\wb_i)$ changes due to update of $\xb$
at iteration $k$. The properties of the conjugate gradient algorithm then do not hold.


The second weakness is that, in most cases (see, e.g.,~\cite{Abubakar_01_TV}.),
the algorithm does not use the exact expression of the gradient
term $\nabla_\xb F$. For all methods based on \eqref{eq:2termes_typique}, we have
\begin{align}
\nabla_\xb F=&~2\lambda\sum_i(\Deltab_i^\dag\Deltab_i\xb-\Deltab_i^\dag\wb_i)+ F_2\nabla_\xb\lambda + \lambda\reg\nabla_\xb \phi(\xb)
\label{eq:gradient_F2_x}
\end{align}
where $\Deltab_i=\diag{\Eb^0_i+\Gv_\cD\wb_i}$ and $\cdot^\dag$ represents the transposed conjugate operation. In general, the term $\nabla_\xb\lambda$ is neglected even though $\lambda$ is a function of $\xb$ according to~\eqref{eq:lambda_CSI}. The
\emph{approximate} gradient is also normalized  by the amplitude of the total field. This latter step is neglected in this paper not to
interfere with our additive regularization term. In practice, we have observed that these approximations prevent the convergence towards a
degenerate minimizer. However, they also prevent the convergence towards a local minimizer, as illustrated in the Section~\ref{results}.
Actually, the solution does not correspond to the stated goal of minimizing $F$.

In the next section, new algorithms designed to avoid the previous pitfalls are proposed.

\section{Proposed algorithms}
\label{sec:propositions} Two new algorithms with improved performance with respect to the original CSI method are now introduced. The first one retains the main
structure of CSI (\ie block-component optimization) while compensating for its deficiencies. The second one is based on a simultaneous
optimization scheme. Both algorithms use the same unique criterion based on \eqref{eq:2termes_typique}.

We first detail the exact form of this criterion. To do so, the expected characteristics of the weight factor are deduced from the optimization
theory and the role of the regularization term is analyzed. Then, we introduce our two optimization strategies.

\subsection{Weight factor}
\label{subsec:facteur_poid} The choice of parameter $\lambda$ is crucial since a bad choice may lead to an inappropriate solution or to a
prohibitive computation time. In this subsection we deduce the expected characteristics of the weight factor from optimization theory and
propose a strategy to set its value. The possible presence of an additional regularization term is deliberately omitted here since it does not
interfere with the weight factor. Regularization issues will be discussed in Subsection~\ref{subsec:reg}.

According to optimization theory~\cite{Nocedal_livre_NO}, if $\lambda\to\infty$ in~\eqref{eq:2termes_typique}, solving~\eqref{eq:F_sum_de_F} is
equivalent to solving
\begin{align}
\label{eq:TMO_opt_sous_cont} \argmin_{\xb,\Wv}F_1\quad\text{under constraint~\eqref{eq:TMO_CS_champs}}
\end{align}
which amounts solving~\eqref{eq:TMO_CS}, or equivalently~\eqref{eq:TMO_non_lin}, in the least-squares sense.  It is therefore quite natural to
consider very large values of $\lambda$ in $F$. Unfortunately, optimization theory~\cite{Nocedal_livre_NO} also states that the minimization
problem~\eqref{eq:F_sum_de_F} becomes ill-conditioned for arbitrary large values of $\lambda$. Practically, this implies that the computation
time rises with $\lambda$.

According to these considerations, $\lambda$ should be set to a value large enough to approximately fulfill the constraint and small enough to
preserve the conditioning of the optimization problem. This trade-off effect and its impact on both solution accuracy and computation time will
be illustrated in the results section.

To our knowledge, no unsupervised method is available to ensure an appropriate choice of $\lambda$. Thus, we rather suggest to turn to a
heuristic tuning. In our experiments, we have observed that an appropriate weight factor for a given contrast is still quite efficient for
"similar" contrasts. We could then imagine that an appropriate $\lambda$ could be set during a training step, involving one or several known
typical contrasts.

We also note that the suggested $\lambda\CSI$ presented in previous section also stems from heuristics since it is not justified by the solution
quality. Hand tuning can then be seen as a generalization of the existing suggestions that gives more flexibility to the user, offering a
trade-off between computation time and accuracy of the solution.

\subsection{Regularization}
\label{subsec:reg}

In this section we analyze the necessity of using a regularization term in criterion $F$.

The continuous-variable microwave tomography problem is intrinsically ill-posed \cite{Isernia_97_NE}. After discretization, it yields an ill-conditioned problem~\cite{Idier_livre_BA} which, in practice, is very sensitive to noise: small perturbations on the measurements cause large variations on the solution.

Regularization, introduced by Tikhonov~\cite{Tikhonov_livre_SI}, is a well known technique to overcome this difficulty. The objective is to
restore the well-posed nature of the problem by incorporating \textit{a priori} information. According to Tikhonov's approach, this is done by
adding a penalty term to the data-adequation component, weighted by a positive parameter $\lambda\reg$ that modulates its relative importance.

The simplest penalty term is the squared norm of the vector of unknowns. Here, we rather penalize the squared norm of first order differences of
$\xb$ between neighboring points in domain $D$, \ie
\begin{equation}
\label{eq:Freg_matrice} F\reg=\phi(\xb)=\| \Dv \xb \|^2
\end{equation}
where each row of $\Dv$ contains only two nonzero values, $1$ and $-1$, in order to properly implement the desired difference.

In a number of existing methods, no penalization is incorporated into the criterion. One could argue that this would lead to degenerate
solutions due to the ill-conditioned nature of the problem. This is not necessarily true if the algorithm is stopped before convergence. Actually, it
has been shown mathematically, in the simpler case of deconvolution-type inverse problems, that stopping a gradient-descent algorithm after $K$
iterations has a similar effect than penalizing the criterion using $F\reg=\| \xb \|^2$, with a weight factor $\lambda\reg$ proportional to
$1/K$~\cite{Lagendijk_88_RI}. We will illustrate this equivalence for the microwave tomography case in the results section.

In the sequel, we shall focus on the penalization approach since it appears as a more explicit way of performing regularization. From the
previous analysis, we choose a criterion based on \eqref{eq:2termes_typique} with \eqref{eq:Freg_matrice} and with $\lambda$ and $\lambda\reg$
constant and set heuristically.

\subsection{Optimization algorithms}
\label{sec:algo_optimi} In this subsection we propose two schemes to minimize the criterion defined previously. The first one, like CSI, is
based on block-component optimization. Meanwhile, it is designed to take full advantage of the conjugate gradient properties. Our second method is based on a
simultaneous update of $\Wv$ and $\xb$. Unlike block-component algorithms, the simplest form of this procedure has the drawback of being
sensitive to the relative scaling of unknowns. Therefore, we propose a technique capable of overcoming this scaling issue.

\subsubsection{Alternated conjugate gradient method}
According to our choices of $\lambda$, $\lambda\reg$ and $F\reg$, the criterion $F$ is quadratic \wrt $\wb_i$ when \xb and $\wb_j$,  $j\neq i$,
are held constant and is also quadratic \wrt $\xb$ when $\Wv$ is held constant.

Indeed, for all $i$, $F$ admits the quadratic expression
\begin{equation}
F(\wb_i)=\wb_i^\dag\Av\wb_i-2\Re(\bb_i^\dag\wb_i)+c_i \label{eq:Fwi}
\end{equation}
as a function of $\wb_i$, where
\begin{align}
\Av=\Gv_\oD^\dag \Gv_\oD+\lambda(\Xv\Gv_\cD-\Iv)^\dag(\Xv\Gv_\cD-\Iv) \label{eq:def_A}
\end{align}
and
\begin{equation}
\bb_i=-\Gv_\oD^\dag\yb_i+\lambda(\Xv\Gv_\cD-\Iv)^\dag\Xv\Eb^0_i \label{eq:def_bi}
\end{equation}
and the quadratic expression
\begin{equation}
F(\xb)=\xb^\dag\Qv\xb-2\Re(\bb^\dag\xb)+c \label{eq:Fx}
\end{equation}
as a function of $\xb$, where
\begin{align}
\label{eq:def_Q} \Qv=\lambda\sum_i\Deltab_i^\dag\Deltab_i+\lambda\reg \Dv^\dag \Dv
\end{align}
and
\begin{equation}
\bb=\lambda\sum_i\Deltab_i^\dag\wb_i\ \ .
\label{eq:def_b}
\end{equation}
Constants $c_i$ and $c$ are not defined here since they are not used in the following analysis.

The linear conjugate gradient algorithm is a gradient-based minimization technique developed to solve linear systems, \ie to minimize quadratic criteria. It has
the advantage of producing computationally low-cost iterations while offering good convergence properties. It is then perfectly suited to
minimize $F$ \wrt each $\wb_i$ and \wrt $\xb$ (the complete form of the linear conjugate gradient algorithm for a complex-valued unknown vector is given in
Appendix \ref{CG}, Table \ref{algo:LCG_general}).

We then propose an algorithm based on two nested iterative procedures. The main loop consists of alternated updates on each sets of unknowns.
Each of these updates is performed by a linear conjugate gradient algorithm. To truly benefit from the efficiency of the conjugate gradient method, we propose to perform
several steps of the algorithm, the first of which not being conjugated with the last one of the previous iteration. This is the main
difference with standard CSI method.

We also use overrelaxation after each update. Convergence of the algorithm is still guaranteed due to the quadratic nature of the subproblems.
Overrelaxation is a well-known \cite{Press92} strategy aimed at limiting the zigzaging effect that will be described in
\ref{subsec:min_groupees}.
The resulting algorithm, which will be referred to as \emph{alternated conjugate gradient for CSI} method, is presented in
Table~\ref{algo:DTC}. Despite nested iterations, it appears that, due to a better use of conjugacy, this algorithm is faster than the standard
CSI method.

\begin{table}[t]
\centering \setlength\fboxsep{2mm}
\begin{spacing}{1.2}
\fbox{\kern-8pt\begin{minipage}{.985\linewidth}
\smallskip
\begin{algorithmic}
\STATE Initialize $\xb=\xb^0$ and $\Wv=\Wv^0$ \STATE $\ell \gets 1$ \REPEAT
\FOR{$i=1,...,M$} \STATE $k \gets 0$ \STATE Initialize $\wb_i$ to $\wb_i^{\ell-1}$ \REPEAT \STATE Perform one iteration of the linear conjugate gradient
algorithm of Table~\ref{algo:LCG_general} to minimize $F$ \wrt $\wb_i$ according to~\eqref{eq:Fwi} \STATE $k \gets k+1$ \UNTIL{Sufficient
decrease of $\norm{\nabla_{w_i}F}^2$} \STATE $\wb_i^\ell \gets \wb_i^{\ell-1} + \theta_\wb(\wb_i-\wb_i^{\ell-1})$ \COMMENT{Overrelaxation} \ENDFOR
\STATE \STATE $k \gets 0$ \STATE Initialize \xb to $\xb^{\ell-1}$ \REPEAT \STATE Perform one iteration of the linear conjugate gradient algorithm of
Table~\ref{algo:LCG_general} to minimize $F$ \wrt $\xb$ according to~\eqref{eq:Fx} \STATE $k \gets k+1$ \UNTIL{Sufficient decrease of
$\norm{\nabla_{\xb}F}^2$} \STATE $\xb^\ell \gets \xb^{\ell-1} + \theta_\xb(\xb-\xb^{\ell-1})$ \COMMENT{Overrelaxation} \STATE $\ell \gets \ell+1$
\UNTIL{Convergence}
\end{algorithmic}
\end{minipage}}
\end{spacing}
\caption{Alternated conjugate gradient for  CSI algorithm} \label{algo:DTC}
\end{table}

As indicated in Table~\ref{algo:DTC}, we propose to perform conjugate gradient steps until a \textit{sufficient decrease} of the gradient is obtained. The
computation cost would be prohibitive if iterations were performed until full convergence (\ie until $\|\nabla F\|=0$). We then suggest to stop
the conjugate gradient algorithm once the initial gradient norm has been reduced by a factor of $T$ \ie to use a \emph{truncated} version of conjugate gradient algorithm.
Typical values of $T$ and overrelaxation coefficients $\theta\in[1,2)$ are $T=10$ or $20$ and $\theta_\wb=\theta_\xb=1.5$.

It should also be underlined that the minimization \wrt $\wb_1\ldots\ \wb_M$ can be performed in parallel, since, according to \eqref{eq:Fwi},
neither $\Av$ nor $\bb_i$ depend on the current values of $\wb_j,\text{ }j \neq i$.

\subsubsection{Simultaneous update algorithm}
\label{subsec:min_groupees}

Both CSI and \emph{alternated conjugate gradient for CSI} algorithms rely on alternated updates of $\Wv$ and $\xb$. However, such schemes are often reported to be inefficient in
classical numerical analysis textbooks such as \cite{Press92} (see Fig.\,10.5.1 therein). More precisely, successive minimization along
coordinate directions suffers from the so-called zigzaging phenomenon: if the minimizer is located within a narrow valley, many small steps are
required to come close to it.

Instead of alternated updates, minimization along conjugate directions is usually recommended. In the microwave tomography context, we are thus driven to apply a
nonlinear conjugate gradient algorithm to the complete set $(\xb,\Wv)$ of unknown quantities. Nonlinear conjugate gradient is an extension of the linear conjugate gradient whose form for
complex-valued unknown vector is detailed in Appendix \ref{CG}, Table \ref{algo:NLCG_general}. The resulting scheme exactly corresponds to
Table~\ref{algo:NLCG_general}, with $\vb=(\xb,\Wv)$ (after arrangement in a column-vector form) and $f=F$. Obviously, we have $\nabla_\vb
f=(\nabla_\xb F,\nabla_{\wb_1} F\ldotsv\nabla_{\wb_M}F)$. A less trivial property is that the optimal step length $\argmin_{\alpha}f(\vb+\alpha
\pb)$ can be computed analytically. It is actually equivalent to finding the zeros of a third order polynomial. The expression of this
polynomial is given in Appendix~\ref{third}.

However, the resulting simultaneous optimization scheme has a drawback compared to alternated optimization schemes: \xb and \Wv correspond to
different physical quantities, which are not measured in the same unit system, and the efficiency of the simultaneous optimization scheme
happens to depend significantly on the chosen units, \ie on the respective scale of \xb and \Wv. Rescaling the variables may result in faster convergence, opur observations indicate that the range for appropriate factors varies widely from experiment to experiment.
Therefore, our approach is to devise a technique that makes the simultaneous optimization scheme insensitive to the relative scales of \xb and
\Wv.  An elegant way to get rid of scaling problems is to use a suitably \emph{preconditioned} version of conjugate gradient algorithm. The general idea
of preconditioning is to resort to the basic conjugate gradient algorithm after a linear invertible change of variables $\vb=\Sv\vb'$, in order to solve an
equivalent but better-conditioned problem~\cite{Nocedal_livre_NO}.  Table~\ref{algo:PCG_general} displays the preconditioned conjugate gradient algorithm
in the case of a complex-valued unknown vector.  The only difference with Table~\ref{algo:NLCG_general} is the presence of the preconditioner
$\Pv=\Sv\Sv^\dag$, which is a positive Hermitian matrix by construction.

\begin{table}[h]
\setlength\fboxsep{2mm}
\begin{spacing}{1.2}
\fbox{\kern-8pt\begin{minipage}{.985\linewidth}
\smallskip
\begin{algorithmic}
\STATE\COMMENT{Minimization of a non quadratic function $f$ \wrt $\vb$} \STATE Initialize \vb \STATE $k \gets 0$ \REPEAT \STATE
$\gb\gets\nabla_\vb f$ \IF{$k=0$} \STATE $\pb\gets-\Pv\gb$ \ELSE \STATE
$\beta\gets\froc{\Re\pth{(\gb-\gb\old)^\dag\Pv\gb}}{(\gb\old)^\dag\Pv\gb\old}$ \STATE $\pb\gets-\Pv\gb+\beta\pb$ \ENDIF \STATE
$\alpha\gets\argmin_{\alpha}f(\vb+\alpha \pb)$ \STATE $\vb\gets\vb+\alpha \pb$ \STATE $\gb\old\gets\gb$ \STATE $k \gets k+1$
\UNTIL{Sufficient decrease of $\norm{\gb}^2$}
\end{algorithmic}
\end{minipage}}
\end{spacing}
\caption{Preconditioned conjugate gradient algorithm in the case of a complex-valued unknown vector} \label{algo:PCG_general}
\end{table}

Here, we use a classical preconditioning strategy in which $\Pv$ is a diagonal matrix formed with the inverse of the diagonal entries of the
Hessian of $F$. Thus, the algorithm becomes independent of a change of units. The proof is rather straightforward and is omitted here.

The main cost of such an algorithm, compared to the unpreconditioned conjugate gradient method, comes from the computation of $\Pv$ which must be done at each iteration. This computation required a matrix-vector multiplication involving a $n\times n$ matrix (as detailed in Appendix \ref{fourth}). In
the results section we will show that despite this additional cost, preconditioned conjugate gradient performs well compared to \emph{alternated conjugate gradient for CSI}.

\begin{changemargin}{0.33cm}{0cm}
\remark{While the conjugate gradient algorithm for a complex-valued unknown vector can always be put into the form of Table~\ref{algo:NLCG_general},
Table~\ref{algo:PCG_general} does not give the general form of preconditioned conjugate gradient in the complex case. The reason is that $\vb=\Sv\vb'$ is not the general
expression for a linear invertible change of variables in the complex case, which would rather read
$[\Re(\vb)\T,\,\Im(\vb)\T]\T=\widetilde{\Sv}[\Re(\vb')\T,\,\Im(\vb')\T]\T$ where $\cdot^\tD$ represents the transpose operation. Therefore, the
complex-valued vectors should be replaced by their equivalent real-valued representation to obtain the general form of preconditioned conjugate gradient. Nonetheless,
preconditioning by a diagonal scaling matrix can always be implemented using complex-valued quantities, so we restrict ourselves to this type of
implementation.}
\end{changemargin}

\section{Results}
\label{results}

We now present some experimental results. First, the pitfalls of the
CSI method and the behavior of the proposed techniques are illustrated
on synthetic data examples. Second, the global performance of the
three studied algorithms (CSI, \emph{alternated conjugate gradient for CSI}
and preconditioned conjugate gradient methods) are compared using
experimental data published in \cite{Belkebir_01_SS}.

Synthetic data were generated with a 2D simulator solving the electric field integral equation using a pulse basis functions and
point-matching scheme. The domain of interest $D$ was of one squared wavelength ($\lambda_0 \times \lambda_0$) of size. We used $M=N$ with
emitters and receivers equally spaced on a circle with radius $\lambda_0/\sqrt{2}$ centered on $D$. Unless otherwise specified, $M=32$ and
$n=32^2=1024$. White Gaussian noise was added to each set of simulated data in order to get a signal-to-noise ratio of 20 dB.

\begin{figure}[t]
\centering \setlength{\tabcolsep}{0pt}
\begin{tabular}{cc}
\begin{tabular}{cc}
\makebox[9pt]{\rotatebox[origin=c]{90}{\footnotesize $\Re\{\xb\}$}}&
\cc{\includegraphics[width=4cm]{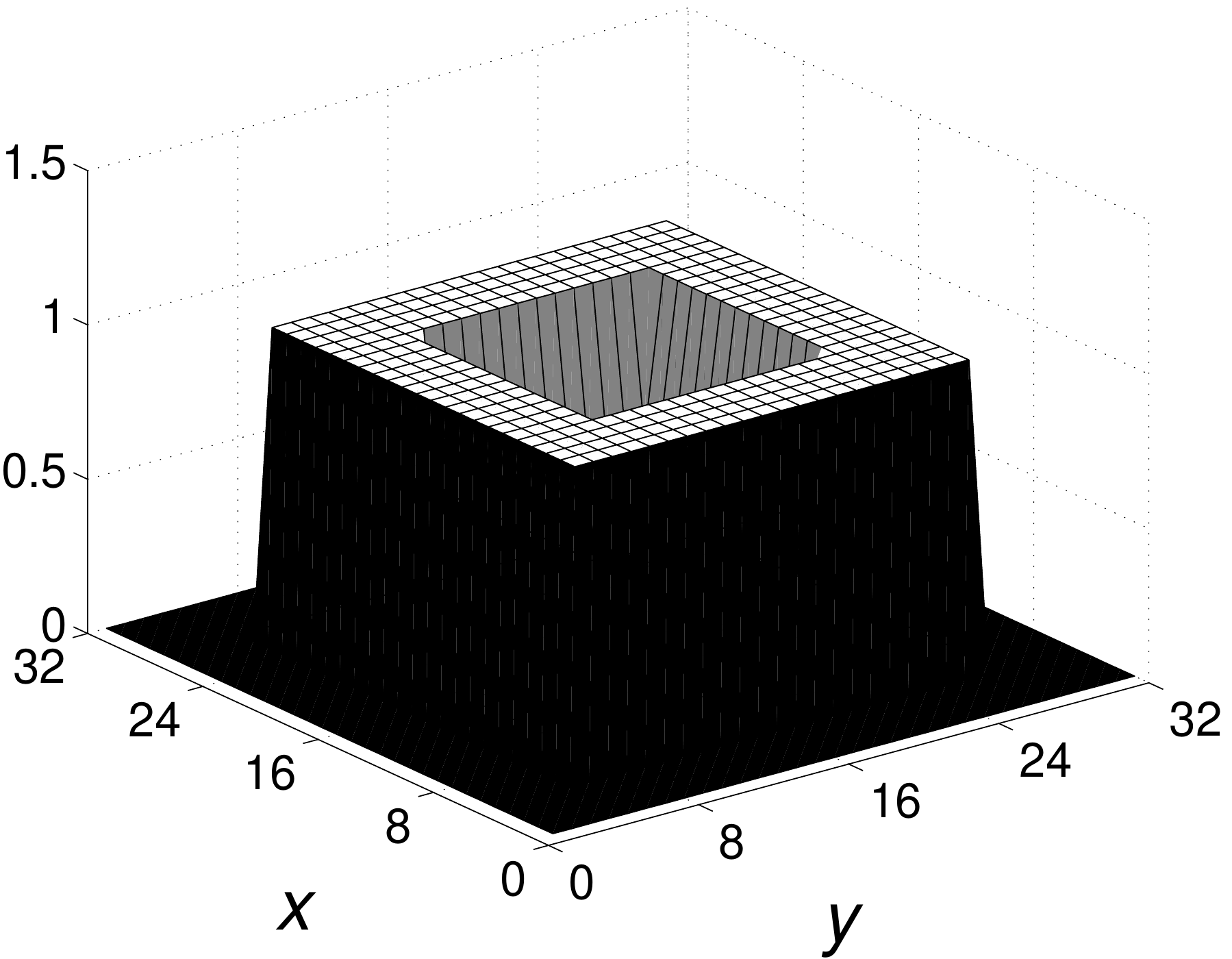}}
\end{tabular}
&
\begin{tabular}{cc}
\makebox[9pt]{\rotatebox[origin=c]{90}{\footnotesize $-\Im\{\xb\}$}}&
\cc{\includegraphics[width=4cm]{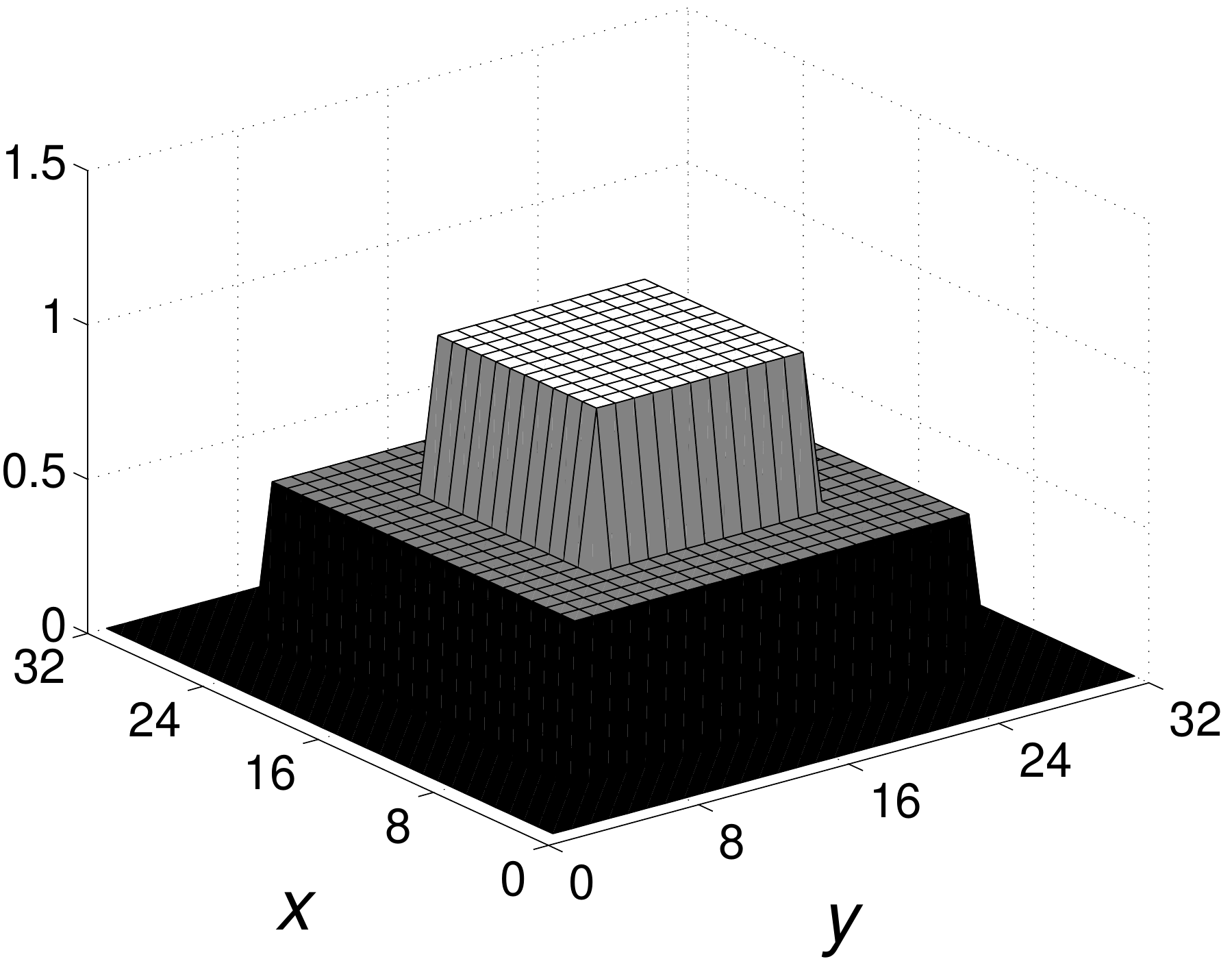}}
\end{tabular}\\
\footnotesize (a) & \footnotesize (b)
\end{tabular}
\caption{(a) Real and (b) minus the imaginary parts of the contrast of the \emph{small square cylinders object}. The \emph{large square cylinder object}  has the same shape but all contrast values are multiplied by a factor of
3. The $x$ and $y$ axes are indexed by the sample number.} \label{fig:OUT_1}
\end{figure}

\begin{figure}[t]
\centering \setlength{\tabcolsep}{0pt}
\begin{tabular}{cc}
\makebox[9pt]{\rotatebox[origin=c]{90}{\footnotesize $\Re\{\xb\}$}}&
\cc{\includegraphics[width=4cm]{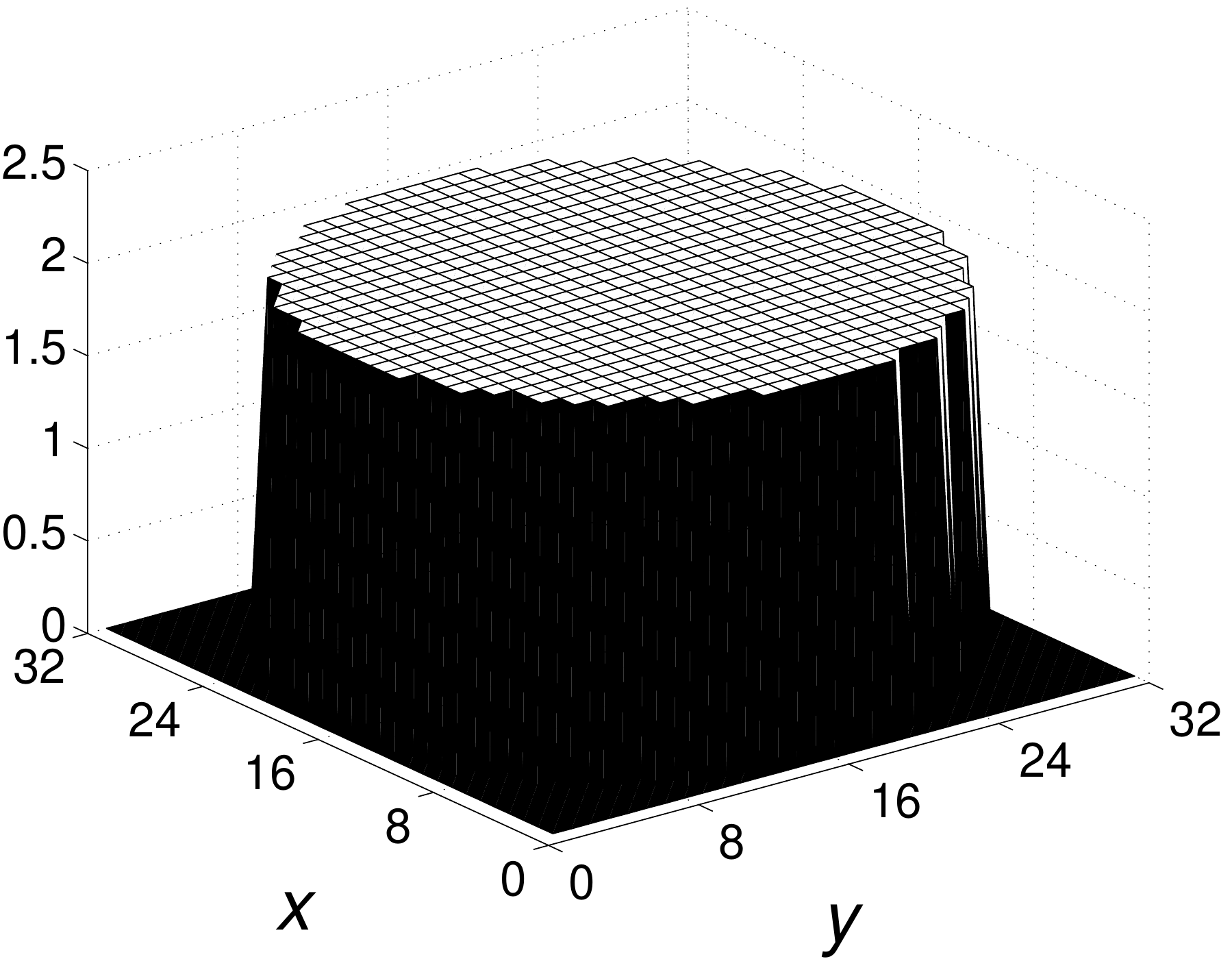}}
\end{tabular}
\caption{Real part of the contrast of the \emph{circular cylinder object}. The imaginary part is zero. The $x$ and $y$ axes are indexed by the sample number.} \label{fig:OUT_3}
\end{figure}

Performance of the algorithms greatly varies according to the shape
and magnitude of the object under test. Tests were thus performed with
three different objects. The first one, shown in Fig.~\ref{fig:OUT_1}, is made of two concentric square cylinders having contrasts of $1-j0.5$, for the outer one, and $0.5-j$, for the inner one. The second object under test has the same shape but all the contrast values are multiplied by a factor 3. These objects will be referred to as the \emph{small square cylinder object} and \emph{large square cylinder object}, respectively. The third object under test, shown in Fig.~\ref{fig:OUT_3}, is made of a single circular cylinder with a radius of $\lambda_0/2$. The contrast of the cylinder is constant and purely real with a value of 2. It will be referred to as the \emph{circular cylinder object}.


Quantitative assessment of the solution quality is performed using a mean square error criterion defined as
\begin{equation}
\label{eq:MSE} \Delta_\xb=\|\x-\x_o\|^2/\|\x_o\|^2
\end{equation}
where $\x_o$ is the actual contrast.

We first illustrate our claims relative to the criterion characteristics and then turn to comparisons related to the optimization process.

\subsection{Criterion characteristics}
\label{subsec:facteur_poid_results} Here we give an example where a local optimization algorithm converges toward a degenerate solution. We
demonstrated the existence of such solutions if $\lambda=\lambda\CSI$ and $\lambda\reg=0$ in \ref{subsec:weak} and Appendix \ref{degenerate}. We
used the \emph{large square cylinders object}  with $N=M=20$ and $n=20^2=400$.

\begin{figure}[t]
\centering \setlength{\tabcolsep}{0pt}
\begin{tabular}{cc}
\makebox[9pt]{\rotatebox[origin=c]{90}{\footnotesize $|\hat{\xb}|$}}&
\cc{\includegraphics[width=6cm]{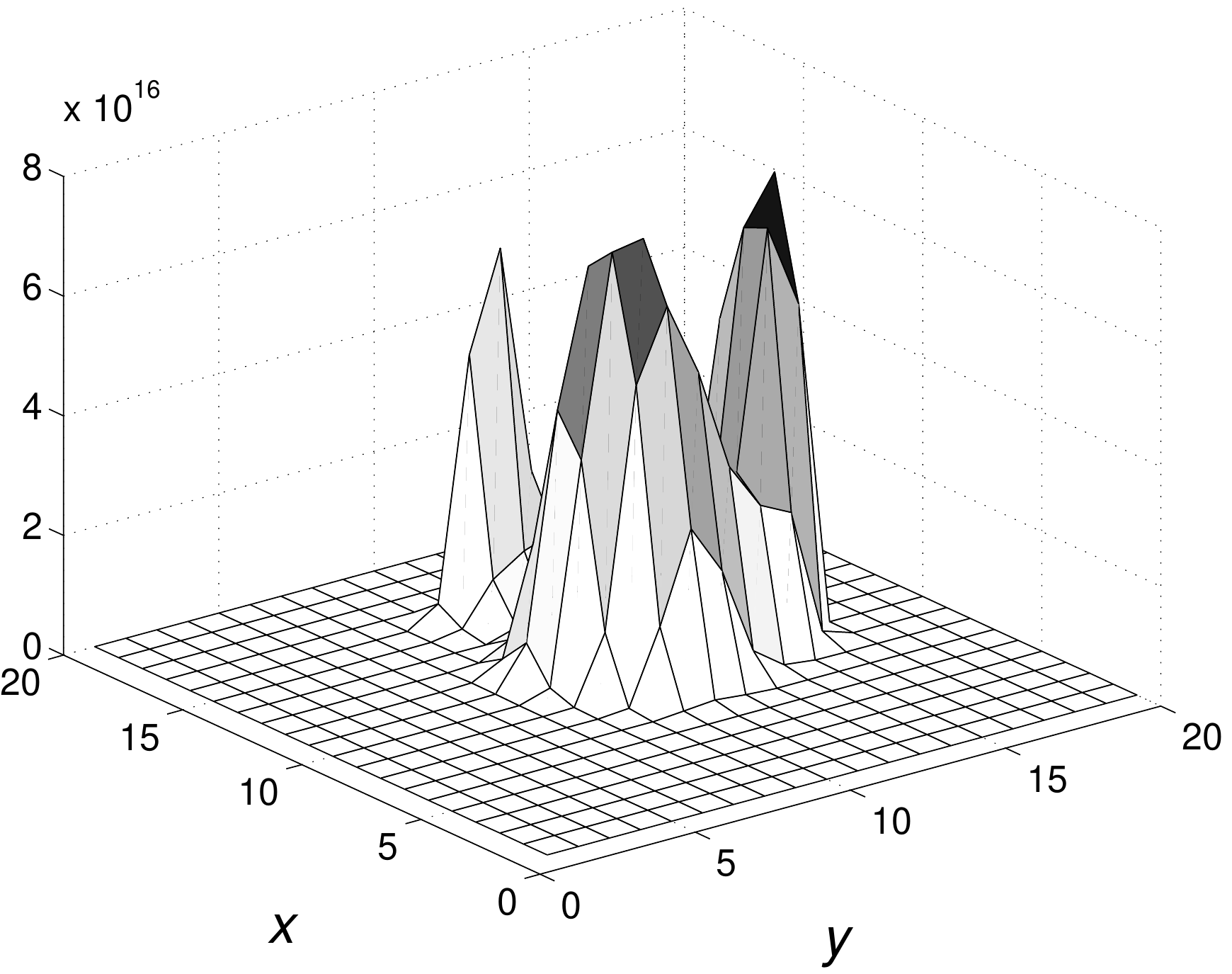}}
\end{tabular}
\caption{Modulus of $\hat{\xb}$ minimizing $F$ with $\lambda=\lambda\CSI$ and $\lambda\reg=0$. \emph{Large square cylinder object}. The $x$ and $y$ axes are indexed by the sample number.} \label{fig:degeneree}
\end{figure}

Fig. \ref{fig:degeneree} presents the modulus of the contrast at the solution while Fig. \ref{fig:total_field_degen} gives the amplitude of
the total field in $D$ for $i=1$ and $5$. As predicted in Appendix \ref{degenerate}, the field vanishes for the pixels $k$ such that
$|\hat{x}_k|\rightarrow\infty$ (for reference, we had $\|\Eb^0_i\|=0.37\ \text{V}/\text{m}$ in the middle of $D$ for each illuminations). This is also true for all
other illuminations. It should be underlined that, in the same conditions, the solution obtained with the \emph{small square cylinder object} and \emph{circular cylinder object} are not degenerate.

\begin{figure}[t]
\centering \setlength{\tabcolsep}{0pt}
\begin{tabular}{cc}
\begin{tabular}{cc}
\makebox[9pt]{\rotatebox[origin=c]{90}{\footnotesize $x$}}&
\cc{\includegraphics[width=4cm]{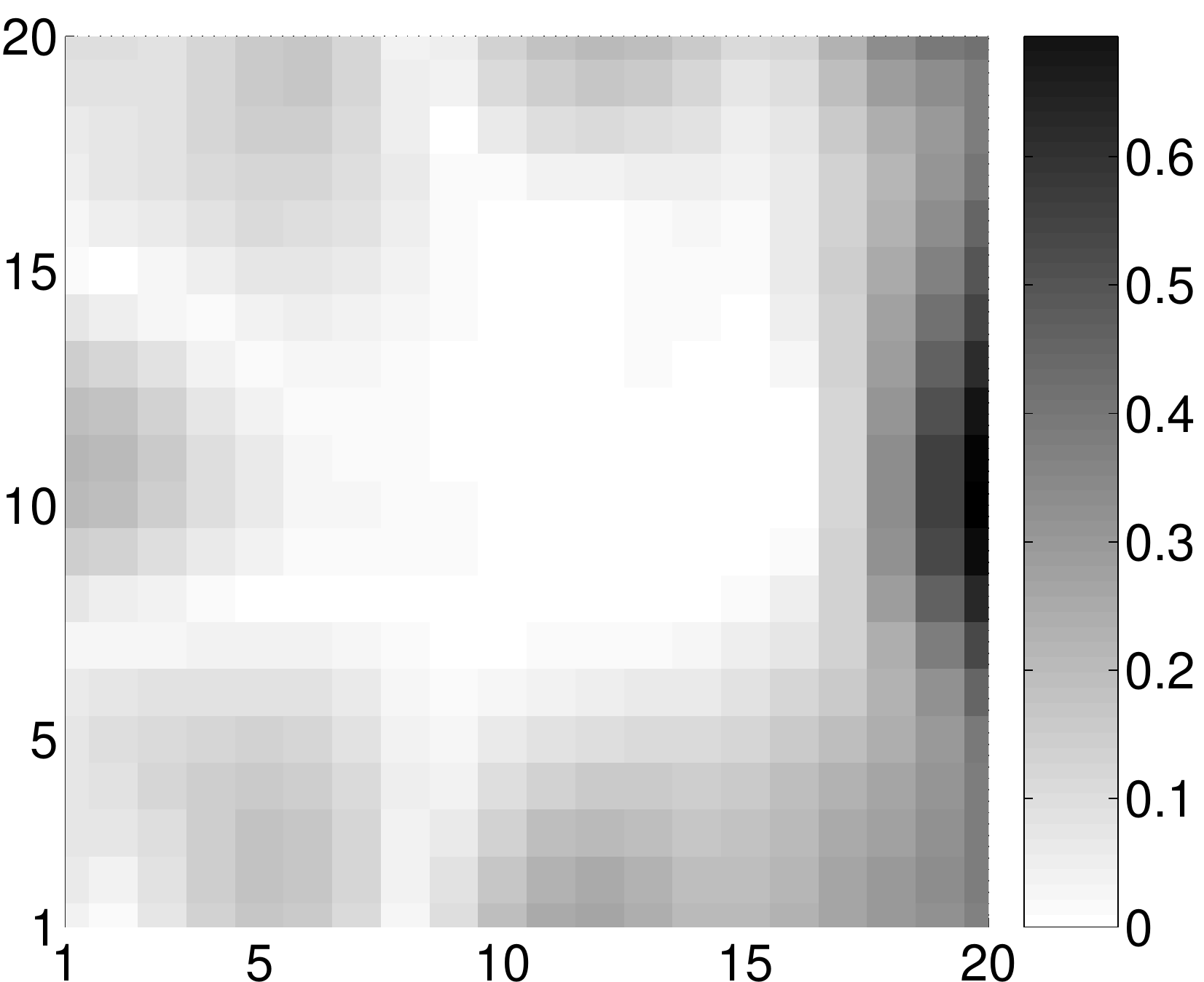}}
\end{tabular}
&
\begin{tabular}{cc}
\makebox[9pt]{\rotatebox[origin=c]{90}{\footnotesize $x$}}&
\cc{\includegraphics[width=4cm]{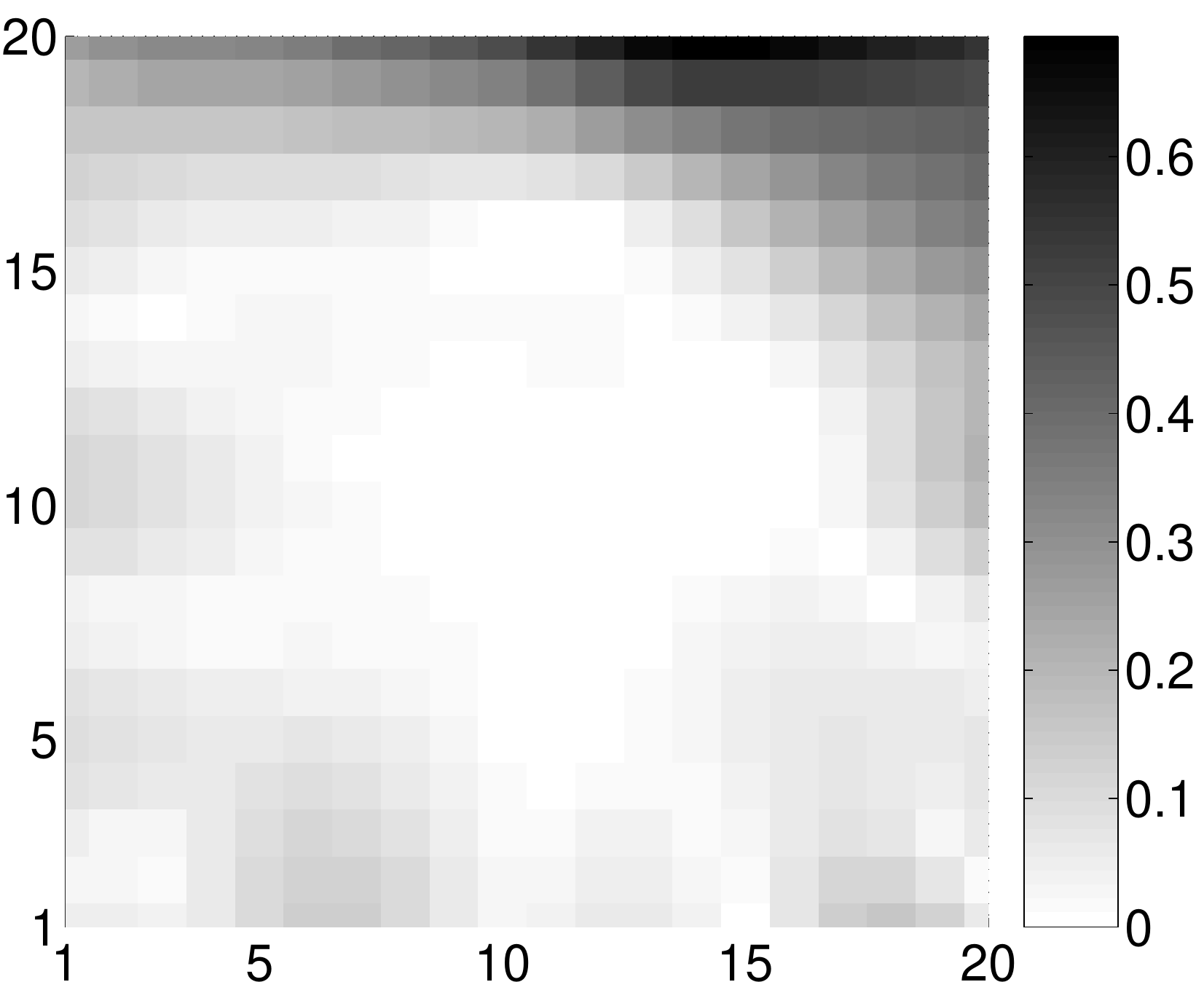}}
\end{tabular}\\
\footnotesize $y$&\footnotesize $y$\\
\footnotesize (a) & \footnotesize (b)
\end{tabular}
\caption{Magnitude of total E-field (in V/m) for the degenerate solution presented in Fig. \ref{fig:degeneree}. Illuminations (a) 1 and (b) 5.
The field vanishes for all pixels $k$ such that $|\hat{x}_k|\rightarrow\infty$. The $x$ and $y$ axes are indexed by the sample number.} \label{fig:total_field_degen}
\end{figure}

\begin{figure}[t]
\centering \setlength{\tabcolsep}{0pt}
\begin{tabular}{ccc}
\makebox[9pt]{\rotatebox[origin=c]{90}{\footnotesize $\Delta_\xb$}}& \cc{\includegraphics[width=6cm]{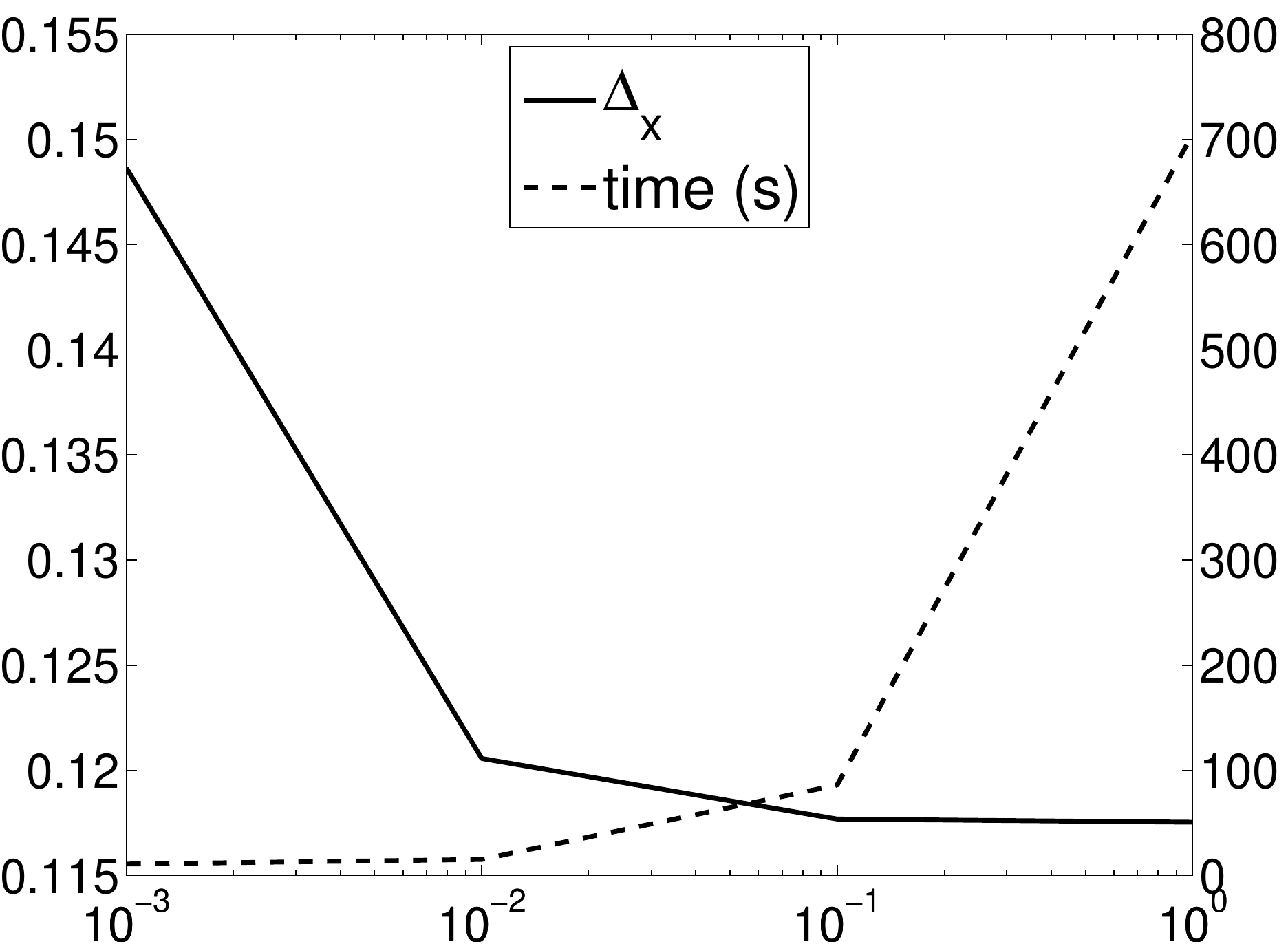}}&
\makebox[9pt]{\rotatebox[origin=c]{90}{\footnotesize Time (s)}}\\
&\footnotesize $\lambda$&
\end{tabular}
\caption{Effect of the weight factor $\lambda$ on the solution quality (mean square error) and on the computation time. The \emph{small square cylinder object} was used with the \emph{alternated conjugate gradient for CSI} method.}
\label{fig:effet_weight_factor}
\end{figure}

According to these results, we suggested, in Subsection \ref{subsec:facteur_poid}, to replace $\lambda\CSI$ by a hand-tuned value of $\lambda$.
In Fig.~\ref{fig:effet_weight_factor} we illustrate the effect of the weight factor on both mean square error and computation time. The \emph{small square cylinder object} and the \emph{alternated conjugate gradient for CSI} method were used.

We remark that, for increasing values of $\lambda$, the computation time rises and quickly becomes prohibitive. On the other hand, we observe
that, for any value of $\lambda$ above a certain threshold, the solution quality remains approximatively unchanged. Moreover, for a given range
of $\lambda$, both solution quality and computation time are acceptable. Within this range, a trade-off can be achieved between computation time
and solution quality.

\begin{figure}[t]
\centering \setlength{\tabcolsep}{0pt}
\begin{tabular}{cc}
\begin{tabular}{cc}
\makebox[9pt]{\rotatebox[origin=c]{90}{\footnotesize $\Re\{\xb\}$}}&
\cc{\includegraphics[width=4cm]{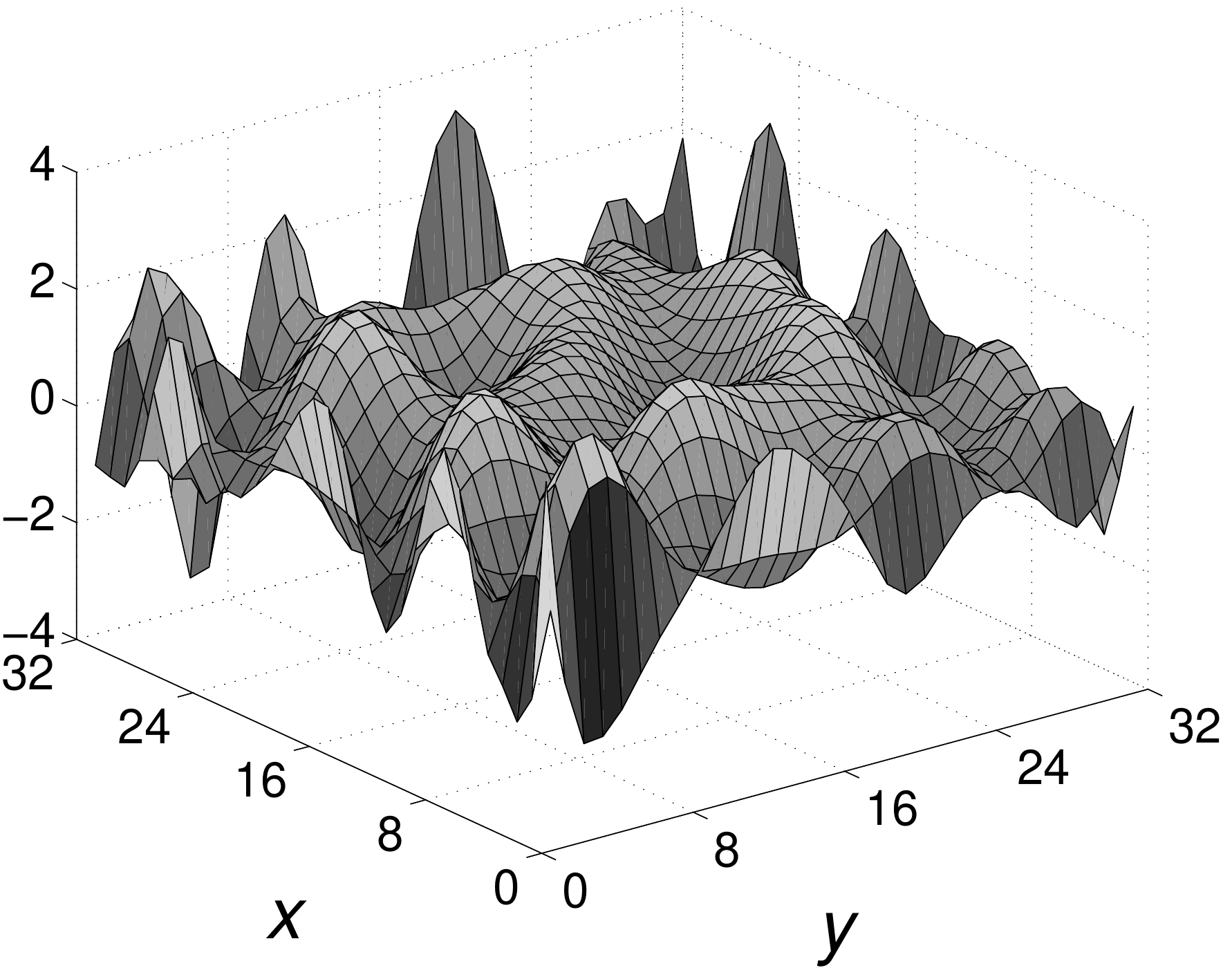}}
\end{tabular}
&
\begin{tabular}{cc}
\makebox[9pt]{\rotatebox[origin=c]{90}{\footnotesize $\Re\{\xb\}$}}&
\cc{\includegraphics[width=4cm]{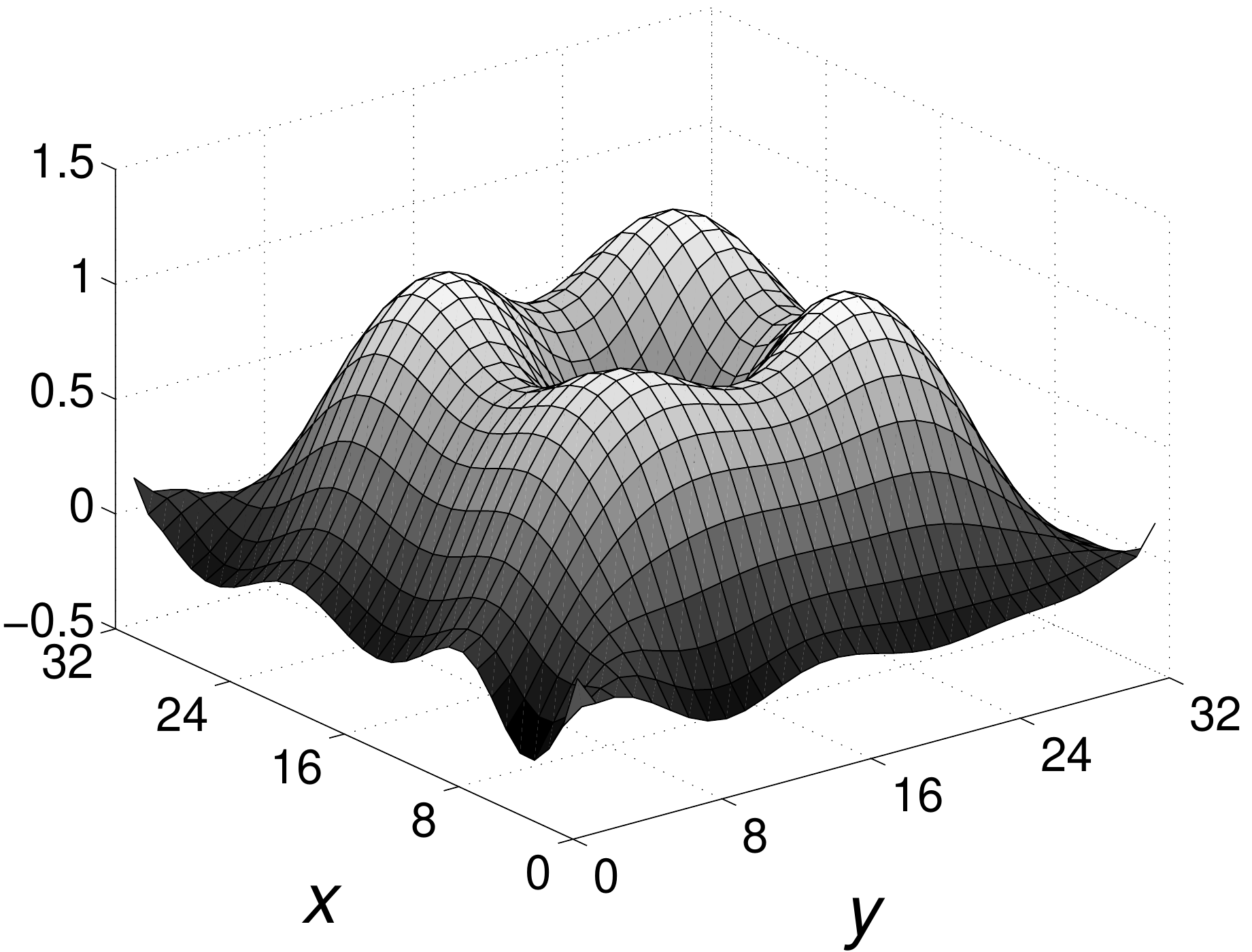}}
\end{tabular}\\
\footnotesize (a) & \footnotesize (b)\\
\multicolumn{2}{c}{
\begin{tabular}{cc}
\makebox[9pt]{\rotatebox[origin=c]{90}{\footnotesize $\Re\{\xb\}$}}&
\cc{\includegraphics[width=4cm]{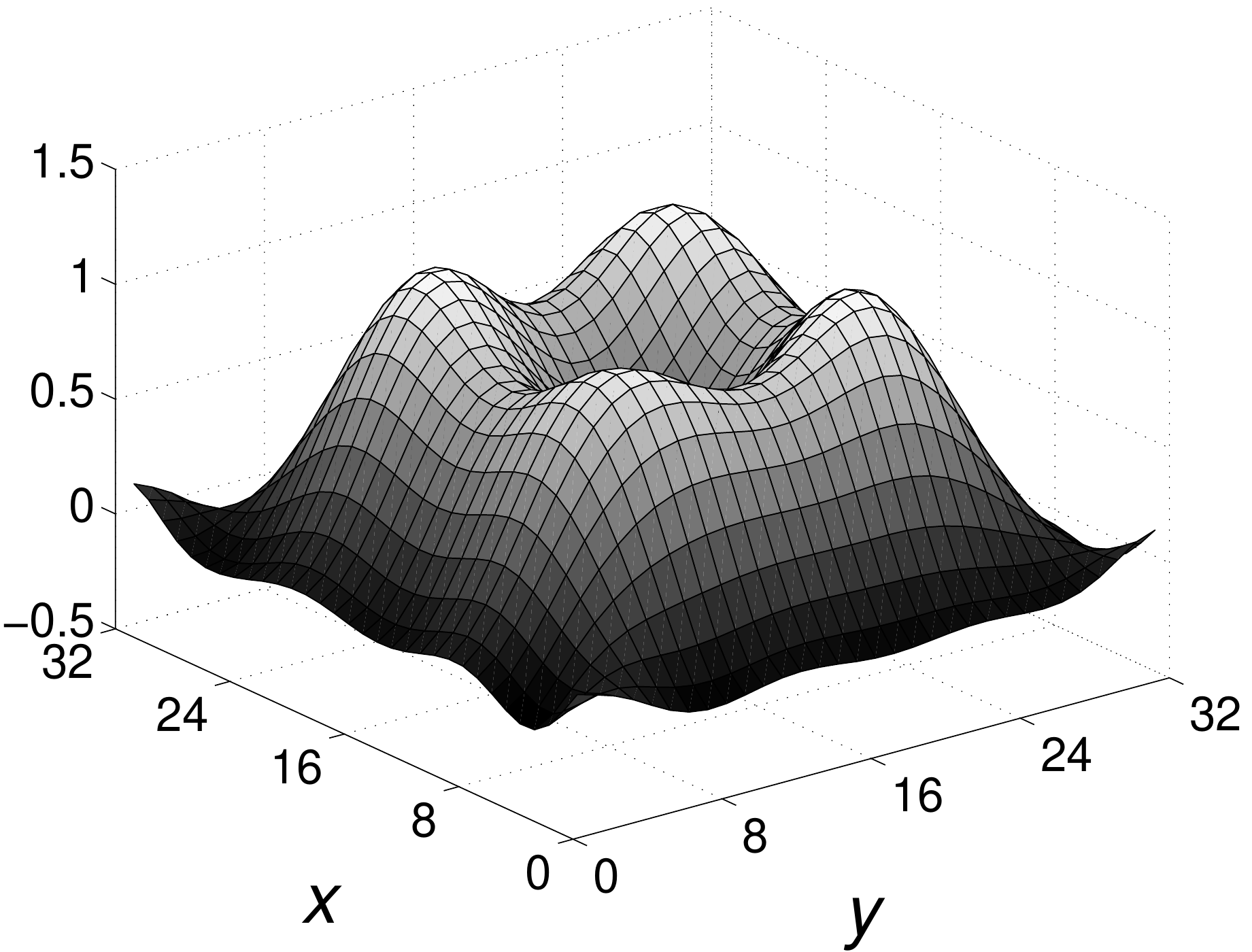}}\end{tabular}}\\
\multicolumn{2}{c}{\footnotesize (c)}
\end{tabular}
\caption{\emph{Small square cylinder object}, real part of the reconstruct contrast: (a) Unregularized criterion, optimization stopped at convergence, (b) unregularized
criterion, optimization stopped before convergence, (c) regularized criterion, optimization stopped at convergence. The $x$ and $y$ axes are indexed by the sample number} \label{fig:effet_reg}
\end{figure}

Finally, we illustrate our assertions of Subsection~\ref{subsec:reg} about regularization. We used the \emph{small square cylinder object} with the \emph{alternated conjugate gradient for CSI} method. Fig.~\ref{fig:effet_reg} (a) presents the solution at convergence (only the real part of the contrast is displayed for the sake of clarity) when an
unregularized criterion is used ($\lambda\reg=0$). We clearly see that the solution is degenerate. In Fig.~\ref{fig:effet_reg} (b) the same
criterion was used but the algorithm was stopped before convergence. Finally, Fig.~\ref{fig:effet_reg} (c) presents the solution obtained at
convergence with a regularized criterion ($\lambda\reg=0.001$). As expected, these two solutions are quite similar.

\subsection{Optimization process}
\label{subsec:Optimization_results}

We now give examples comparing the CSI optimization scheme to the proposed methods.

\begin{figure}[t]
\centering \setlength{\tabcolsep}{0pt}
\begin{tabular}{cc}
\makebox[9pt]{\rotatebox[origin=c]{90}{\footnotesize $\|\nabla_xF\|$}}&
\cc{\includegraphics[width=6cm]{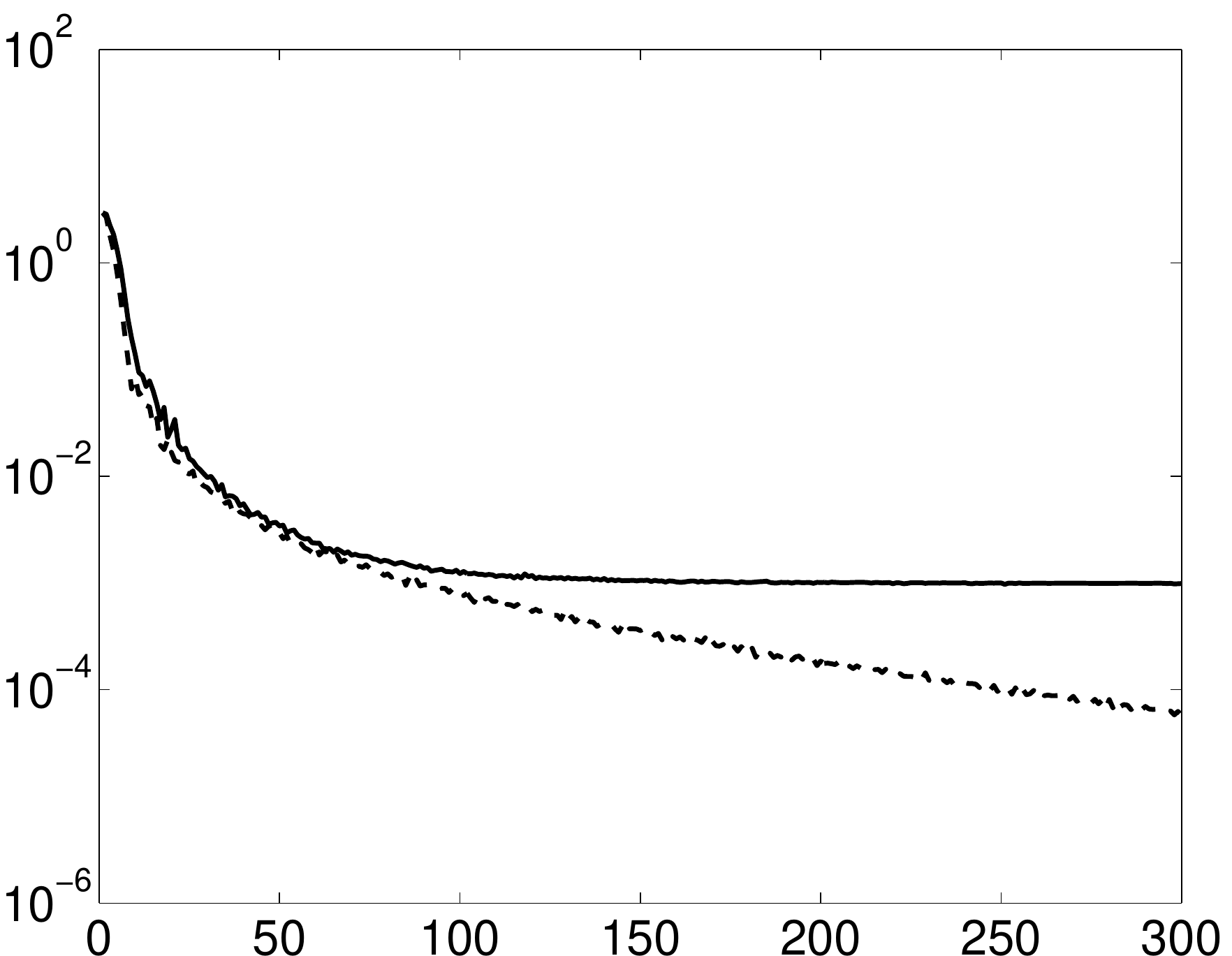}}\\
&\footnotesize Iteration number
\end{tabular}
\caption{Evolution of norm of the gradient \wrt $\xb$ for the CSI method when (-) a gradient approximation is used to calculate update
directions and when (-~-) no approximation is used. Remark: the plotted quantity is the exact value of the norm of the gradient, not its
approximation.} \label{fig:convergence_CSI}
\end{figure}

In Subsection \ref{subsubsec:CSI_algo_weak} we stated without proof that the approximations on which standard CSI relies prevent the convergence
toward a local minimizer of $F$. This point is illustrated in Fig. \ref{fig:convergence_CSI}, which depicts the evolution of norm of the gradient as a function of the number of iterations. With the CSI approximations~\cite{Abubakar_01_TV} (solid lines) the gradient norm goes toward a nonzero value, thereby indicating that the convergence point is not a local minimizer. When no approximations are made (dashed lines), the gradient norm decreases toward zero as expected.


Another comment in Subsection \ref{subsubsec:CSI_algo_weak} was related to the non-optimal exploitation of the conjugate gradient algorithm in the standard CSI
optimization scheme. Here we validate that the \emph{alternated conjugate gradient for CSI} scheme is actually more efficient than the original CSI: both methods were tested using the same
criterion. Parameters $\lambda$ and $\lambda\reg$ were set by hand to $0.01$ and $0.001$, respectively. Tests were performed with the \emph{small square cylinder object} and on the \emph{circular cylinder object}.

\begin{figure}[t]
\centering \setlength{\tabcolsep}{0pt}
\begin{tabular}{cc}
\begin{tabular}{cc}
\makebox[9pt]{\rotatebox[origin=c]{90}{\footnotesize $F$}}&
\cc{\includegraphics[width=4cm]{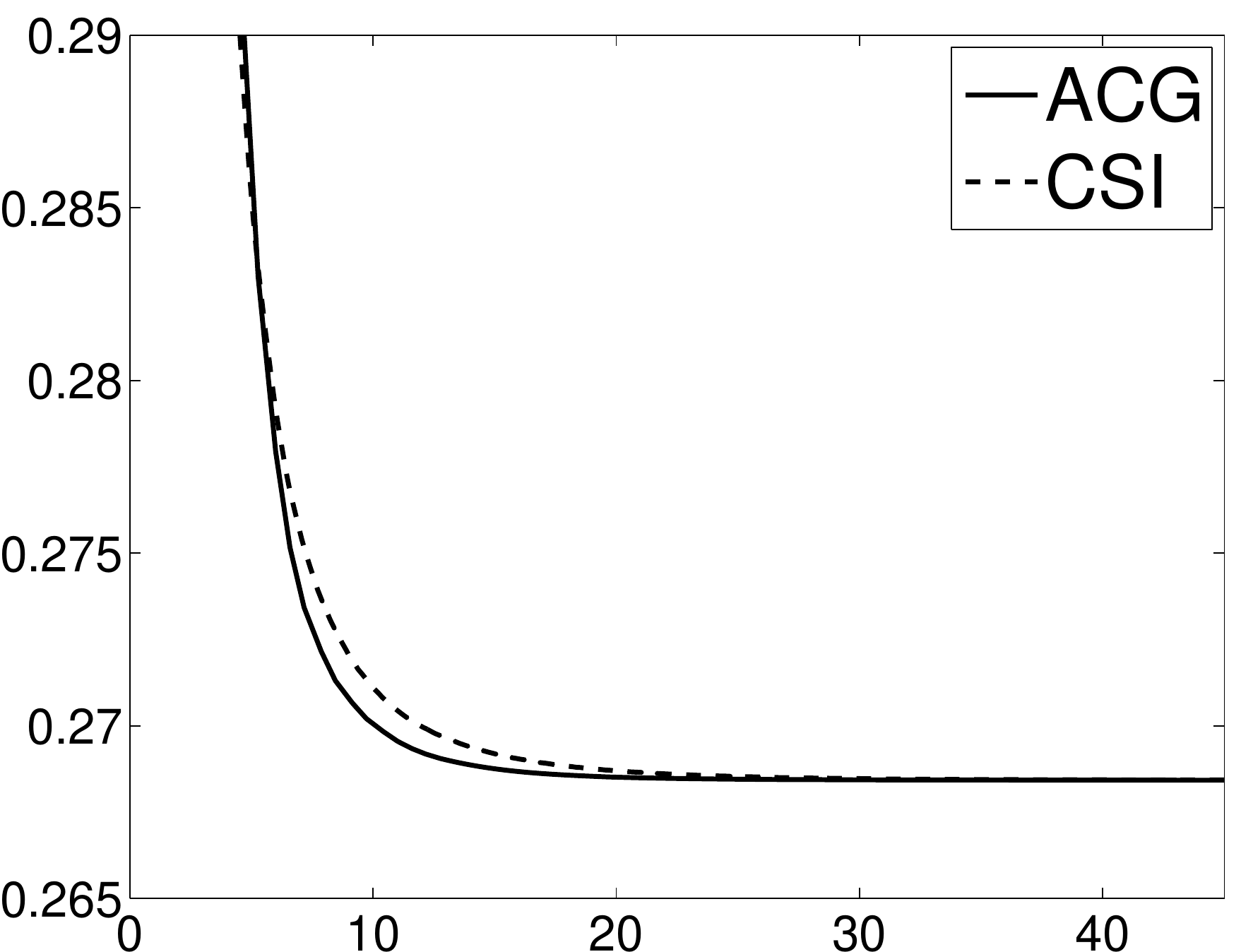}}\\
&\footnotesize Time (s)
\end{tabular}
&
\begin{tabular}{cc}
\makebox[9pt]{\rotatebox[origin=c]{90}{\footnotesize $F$}}&
\cc{\includegraphics[width=4cm]{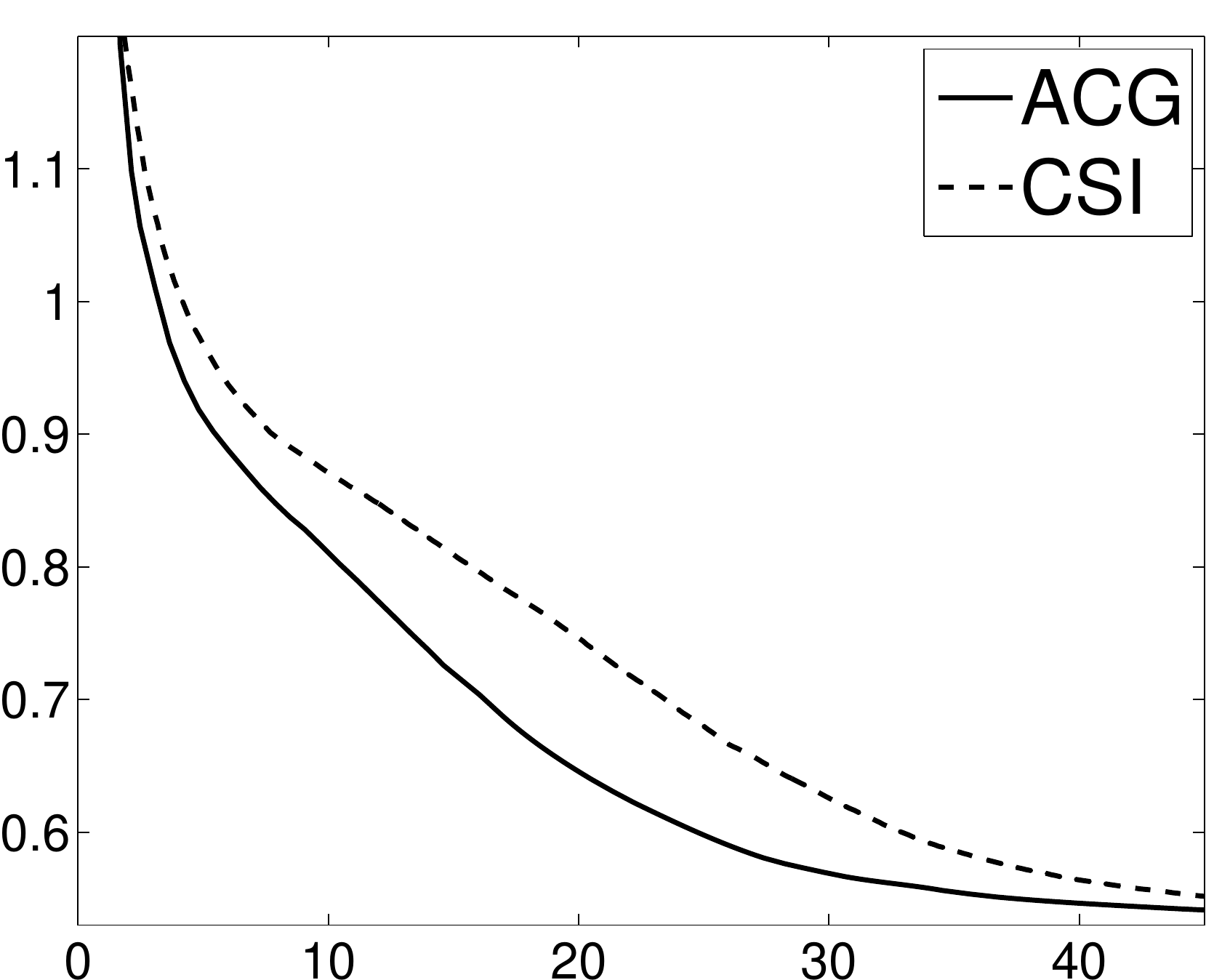}}\\
&\footnotesize Time (s)
\end{tabular}\\
\footnotesize (a) & \footnotesize (b)
\end{tabular}
\caption{Evolution of $F$ as a function of time for CSI and \emph{alternated conjugate gradient for CSI} (ACG) schemes: (a) \emph{small square cylinder object}, (b) \emph{circular cylinder object}. The same criterion, with $\lambda=0.01$ and
$\lambda\reg=0.001$, was used for both methods.} \label{fig:CSI_vs_DTC}
\end{figure}

Figs. \ref{fig:CSI_vs_DTC} (a) and (b) present the evolution of criterion $F$ as a function of time for the \emph{small square cylinder object} and the \emph{circular cylinder object}, respectively. In both
cases, the \emph{alternated conjugate gradient for CSI} algorithm is faster, suggesting that a better use of conjugacy pays off. Experience shows that the gain in computation time can
be as high as $20\%$ depending on the object under test and on the chosen stopping rule. Obviously, it cannot be proved that the \emph{alternated conjugate gradient for CSI} scheme is always
faster than CSI. Nevertheless, \emph{alternated conjugate gradient for CSI} outperformed CSI in all test cases.

In Subsection \ref{subsec:min_groupees} we proposed algorithms performing simultaneous updates of the unknowns. We illustrate how these
simultaneous conjugate gradient and preconditioned conjugate gradient algorithms behave and we compare them to the \emph{alternated conjugate gradient for CSI} method.

\begin{figure}[t]
\centering \setlength{\tabcolsep}{0pt}
\begin{tabular}{cc}
\begin{tabular}{cc}
\makebox[9pt]{\rotatebox[origin=c]{90}{\footnotesize $F$}}&
\cc{\includegraphics[width=4cm]{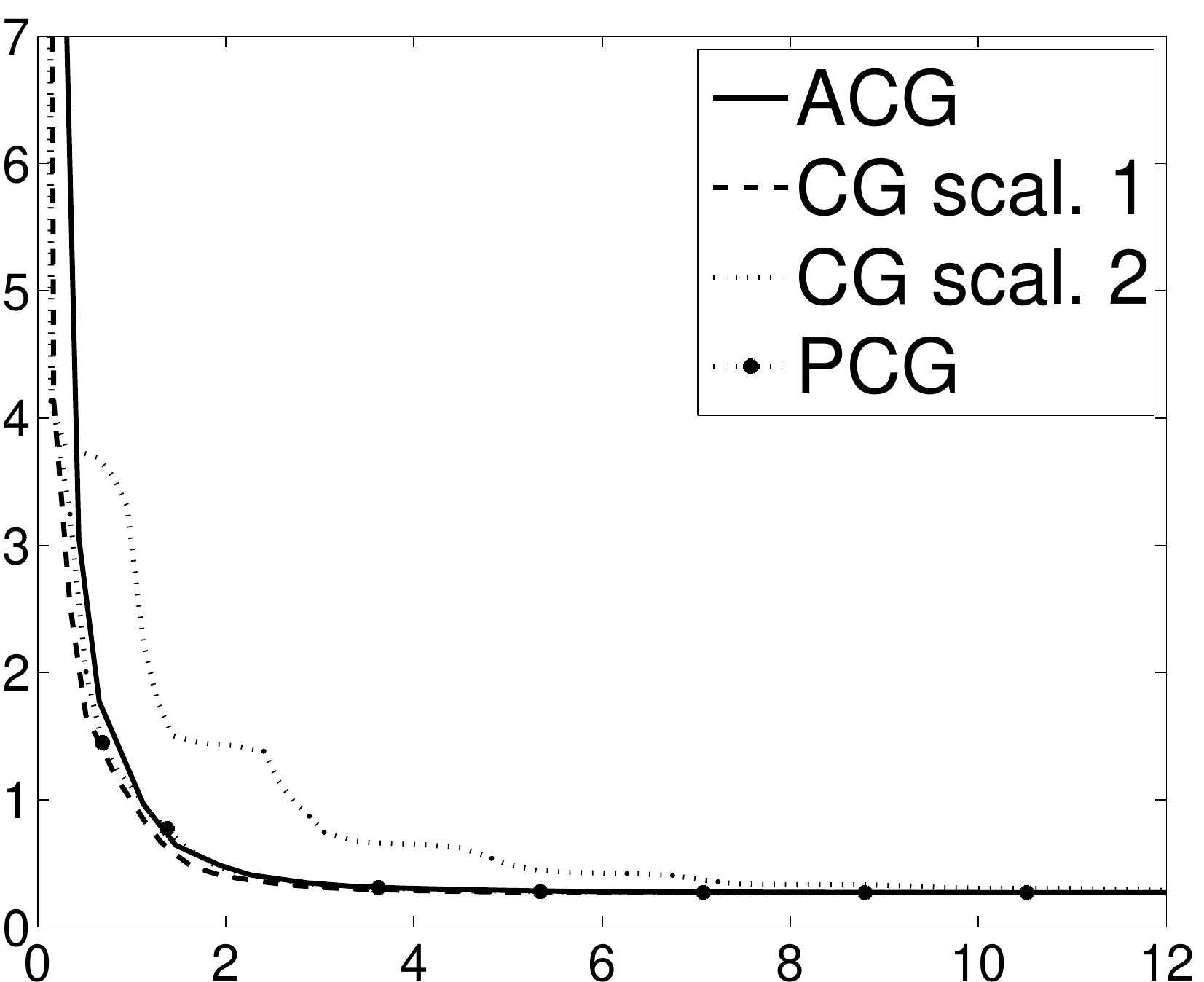}}\\
&\footnotesize Time (s)
\end{tabular}
&
\begin{tabular}{cc}
\makebox[9pt]{\rotatebox[origin=c]{90}{\footnotesize $F$}}&
\cc{\includegraphics[width=4cm]{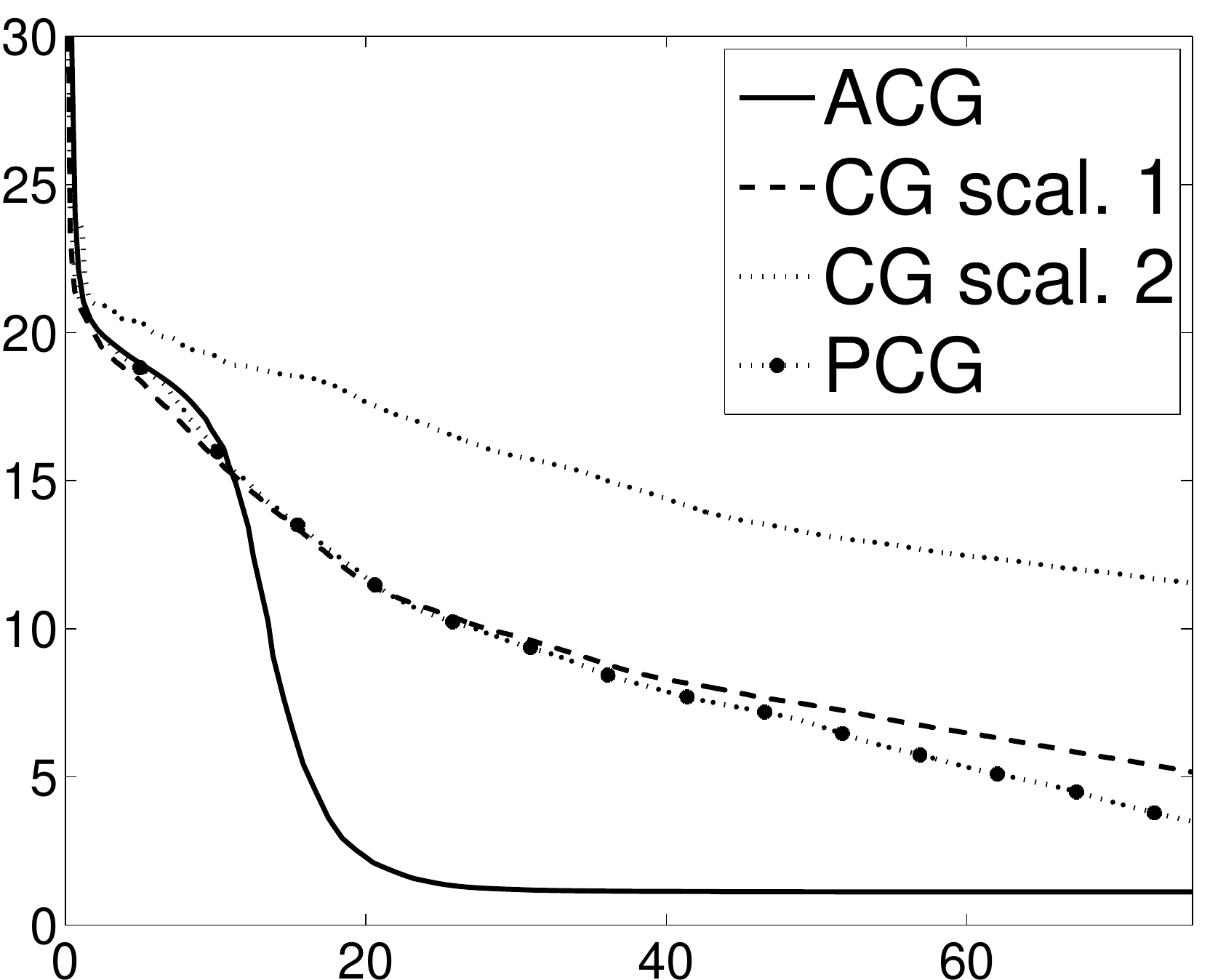}}\\
&\footnotesize Time (s)
\end{tabular}\\
\footnotesize (a) & \footnotesize (b)
\end{tabular}
\caption{Evolution of $F$ as a function of time for: \emph{alternated conjugate gradient for CSI} (ACG) algorithm with scaling 1, preconditioned conjugate gradient (PCG) algorithm with scaling 1, unpreconditioned conjugate gradient (CG) algorithm
with scaling 1 and unpreconditioned conjugate gradient algorithm with scaling 2. For preconditioned conjugate gradient and \emph{alternated conjugate gradient for CSI} algorithms, the results are identical for scalings 1 and 2.
(a) \emph{small square cylinder object}, (b) \emph{circular cylinder object}} \label{fig:simul_vs_DTC}
\end{figure}

In Fig. \ref{fig:simul_vs_DTC}, the evolution of the criterion is
presented for unpreconditioned conjugate gradient, preconditioned
conjugate gradient and \emph{alternated conjugate gradient for CSI} methods. Two different scales
were used. The results in Fig. \ref{fig:simul_vs_DTC} (a) were obtained with the \emph{small square cylinder object}  while those in Fig. \ref{fig:simul_vs_DTC} (b) were obtained
with the \emph{circular cylinder object}. Scaling~2  differs from
scaling~1  by
the fact that a factor of one tenth was applied to the currents.

The scaling sensitivity of the unpreconditioned conjugate gradient algorithm appears clearly in the results. The convergence speed  of the algorithm exhibits large
variations when units are changed. This is not the case for the other two methods. Moreover, according to our experience, the preconditioned conjugate gradient
method always provides results at least as good as those produced by the unpreconditioned conjugate gradient technique, and should therefore be preferred.

However, comparison of simultaneous preconditioned conjugate gradient and block-component optimization is rather inconclusive, as the nature of the object under test seems to have a significant impact. Indeed, while the simultaneous scheme is a little bit faster for the \emph{small square cylinder object}, its convergence speed is not even competitive for
the \emph{circular cylinder object}. Those variations prevent us from systematically favoring one type of algorithm over the other. More
details are given on this subject in the next subsection.

\subsection{Experimental data}
\begin{figure}[t]
\centering \setlength{\tabcolsep}{0pt}
\begin{tabular}{cc}
\begin{tabular}{cc}
\makebox[9pt]{\rotatebox[origin=c]{90}{\footnotesize $\Re\{\xb\}$}}&
\cc{\includegraphics[width=4cm]{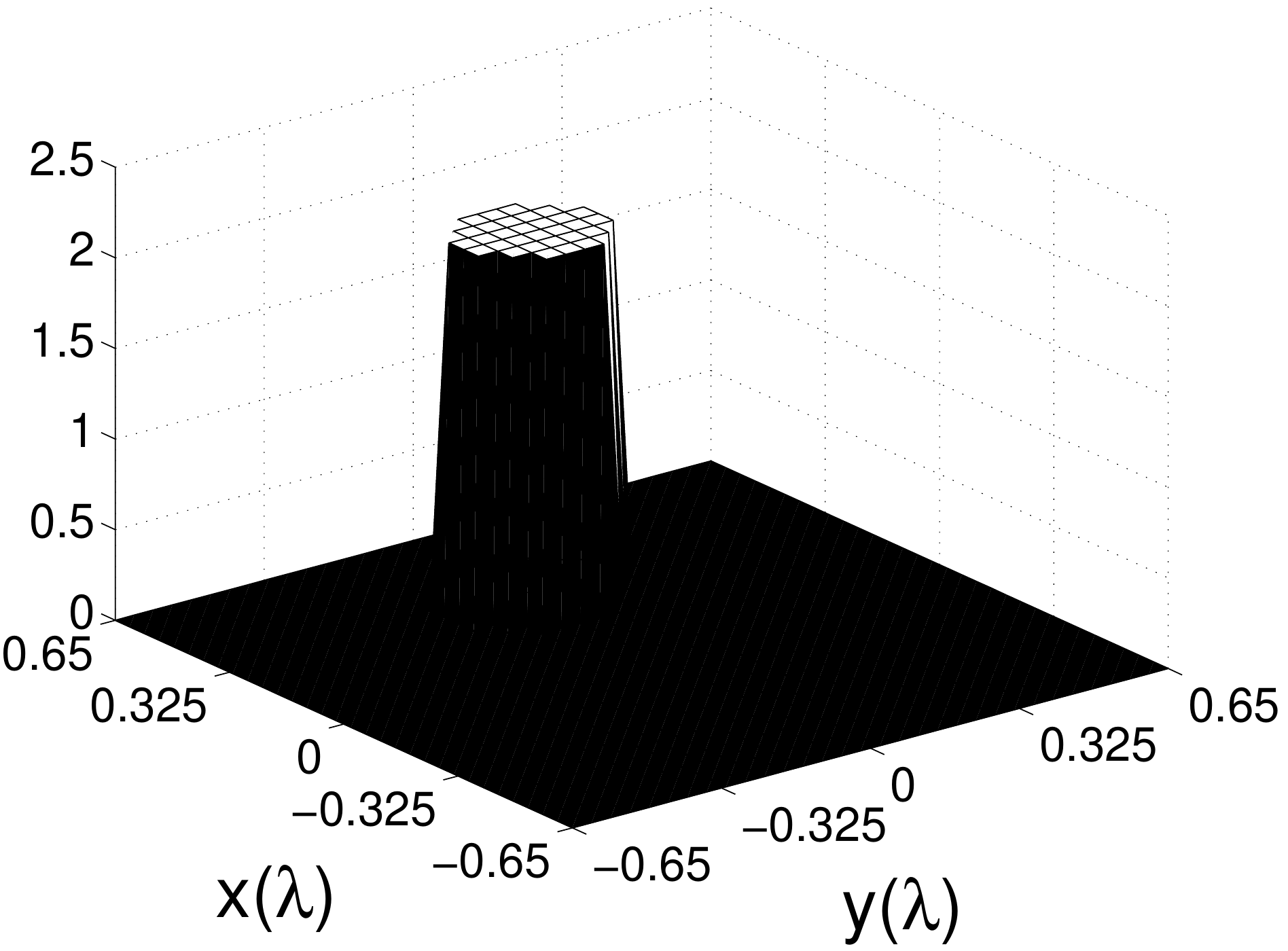}}
\end{tabular}
&
\begin{tabular}{cc}
\makebox[9pt]{\rotatebox[origin=c]{90}{\footnotesize $\Re\{\xb\}$}}&
\cc{\includegraphics[width=4cm]{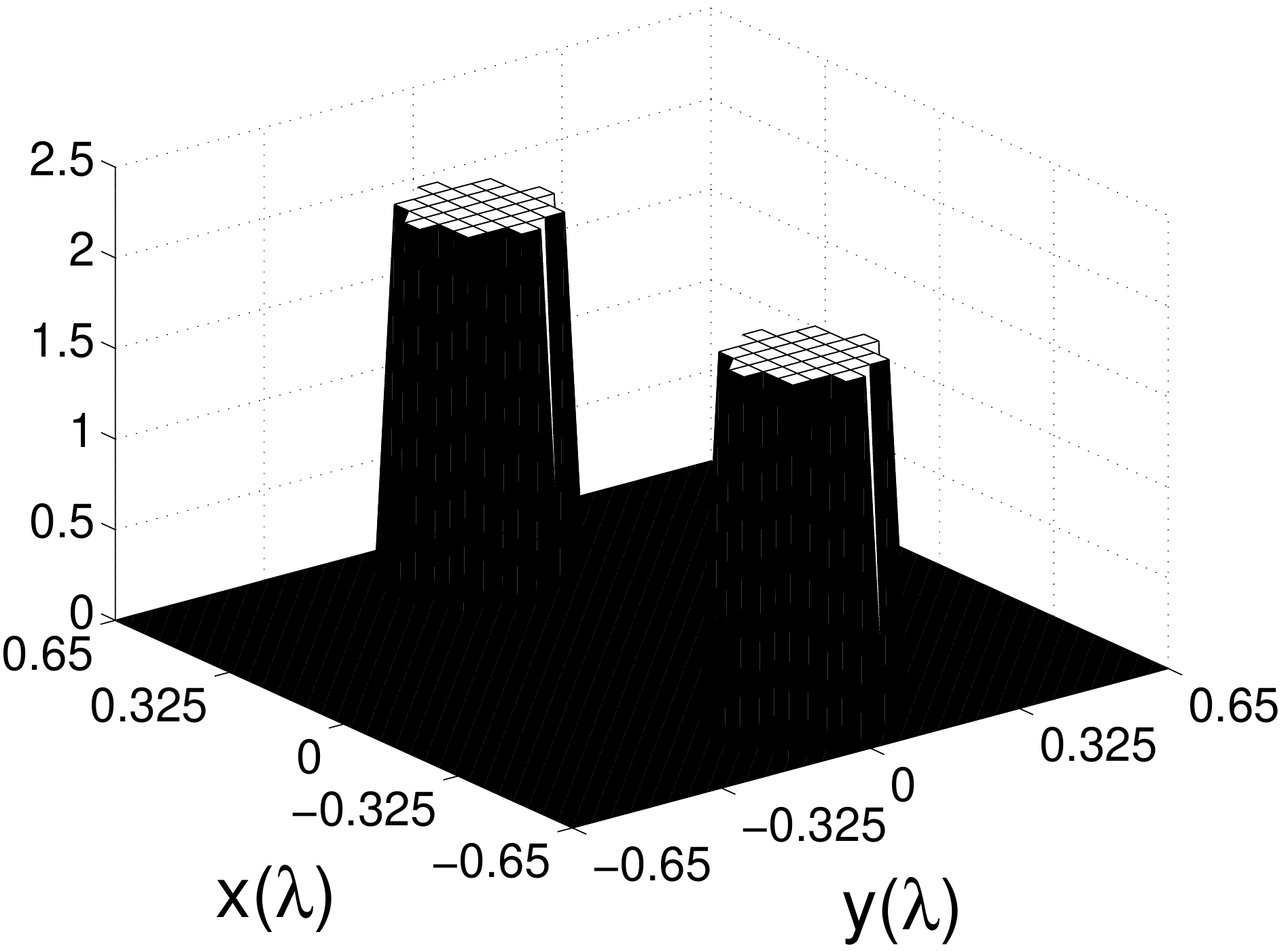}}
\end{tabular}\\
\footnotesize (a) & \footnotesize (b)
\end{tabular}
\caption{Real parts of the contrasts of objects under test tested in \cite{Belkebir_01_SS}. (a) \emph{One cylinder object} (b) \emph{two cylinder object}. Both objects under test have a purely real
contrast.} \label{fig:fresnel_OUT}
\end{figure}

We now compare the studied algorithms using the experimental data published in~\cite{Belkebir_01_SS} (the $3$ GHz dataset was used). In this paper, a quasi 2-D setup is used
to perform measurements over two different objects under test, presented in Fig.~\ref{fig:fresnel_OUT}~(a) and (b). They will be referred to  as the \emph{one cylinder}
and \emph{two cylinder} cases, respectively. Their contrasts are purely real. For more detail on the setup, see \cite{Belkebir_01_SS}.

In our tests, parameter $\lambda$ was set as follows: For CSI, we always used $\lambda\CSI$ as defined in \eqref{eq:lambda_CSI}. For \emph{alternated conjugate gradient for CSI} and
preconditioned conjugate gradient, two values were used: the final value of $\lambda\CSI$ (test~1); a heuristic value offering a good trade-off between solution quality and
computation time (test~2). Finally, we used the same regularization factor $\lambda\reg$ for all tests and all three methods.



\begin{figure}[t]
\centering \setlength{\tabcolsep}{0pt}
\begin{tabular}{cc}
\begin{tabular}{cc}
\makebox[9pt]{\rotatebox[origin=c]{90}{\footnotesize $\Delta_\xb$}}&
\cc{\includegraphics[width=4cm]{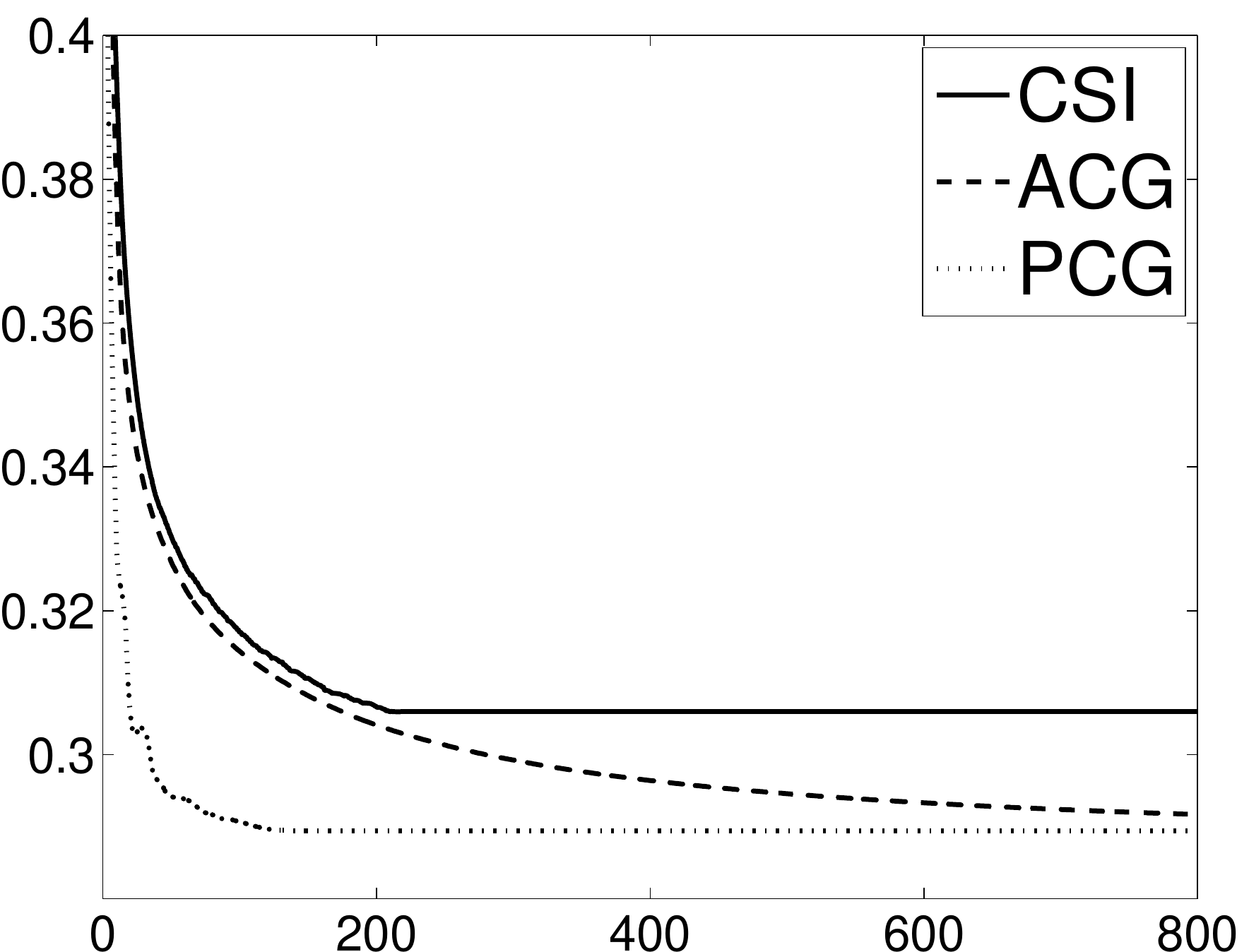}}\\
&\footnotesize Time (s)
\end{tabular}
&
\begin{tabular}{cc}
\makebox[9pt]{\rotatebox[origin=c]{90}{\footnotesize $\Delta_\xb$}}&
\cc{\includegraphics[width=4cm]{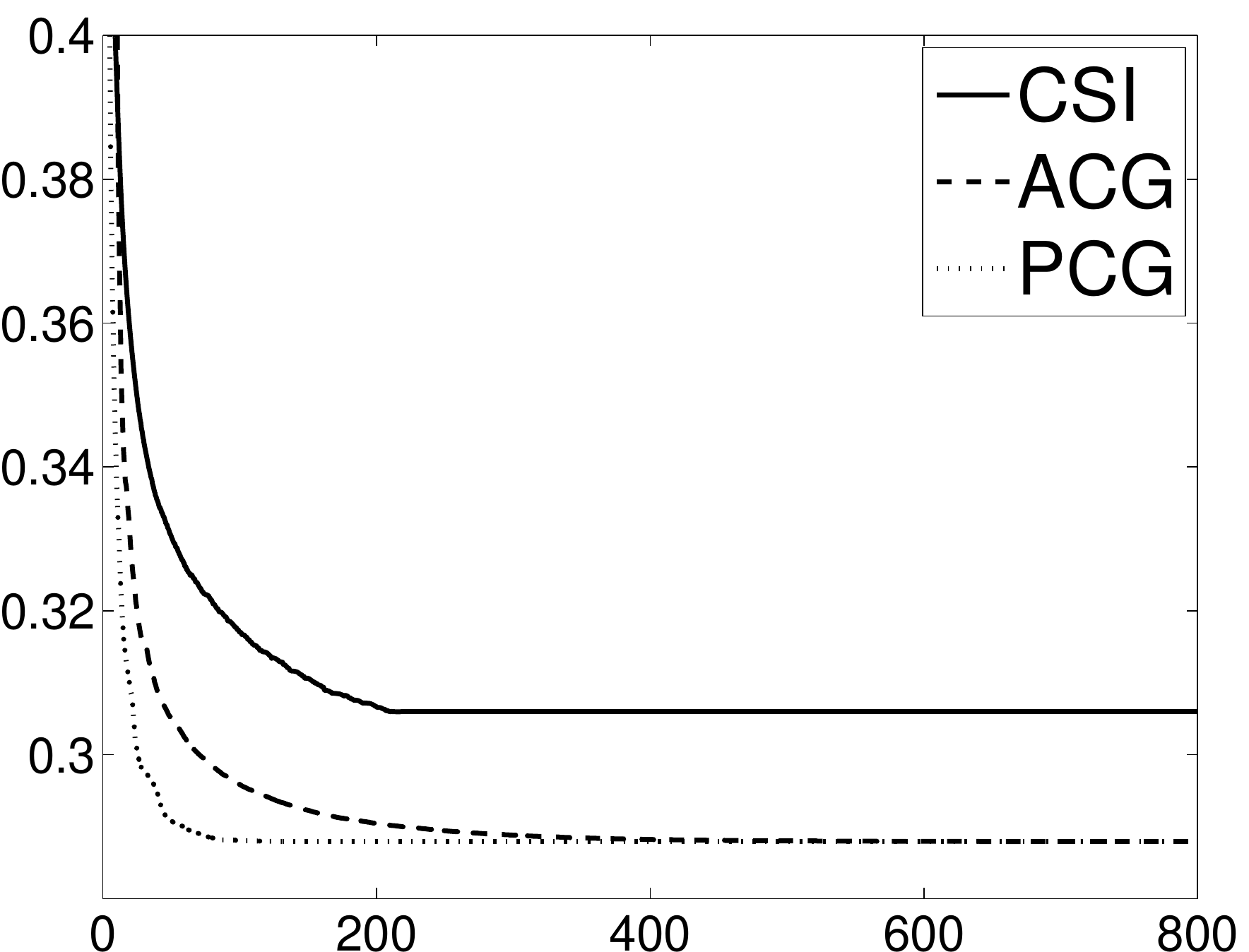}}\\
&\footnotesize Time (s)
\end{tabular}\\
\footnotesize (a) & \footnotesize (b)
\end{tabular}
\caption{Evolution of $\Delta_\xb$ in function of time for CSI algorithm, \emph{alternated conjugate gradient for CSI} (ACG) method and preconditioned conjugate gradient (PCG) algorithm. $\lambda =\lambda\CSI$ for the CSI
method. For \emph{alternated conjugate gradient for CSI} and preconditioned conjugate gradient methods: (a) $\lambda$ equal the final value of $\lambda\CSI$ , (b) $\lambda$ set heuristically. Data from~\cite{Belkebir_01_SS}.}
\label{fig:fresnel_results}
\end{figure}

Fig. \ref{fig:fresnel_results} presents the evolution of $\Delta_\xb$ as a function of time (the evolution of $F$ cannot be used  here for
comparison purpose since the criteria are not the same for all algorithms). Figs. \ref{fig:fresnel_results} (a) and \ref{fig:fresnel_results}
(b) present the results for test~1 and 2, respectively.


The results of the first test show that \emph{alternated conjugate gradient for CSI} and preconditioned conjugate gradient provide better solutions than CSI. This  seems to contradict the fact that the values of $\lambda$ were selected so as to obtain identical criterion values at the solution point.
Actually, this difference can be attributed to the fact that the CSI algorithm does not converge toward a local minimum of the criterion, as
illustrated in Subsection \ref{subsec:Optimization_results}.

Moreover, the \emph{alternated conjugate gradient for CSI} method is faster than CSI in the first test. This once again suggests that the better use of conjugacy made in the \emph{alternated conjugate gradient for CSI}
method pays off.

The second test also reveals that the hand tuning of $\lambda$ yields a significant increase of the speed of \emph{alternated conjugate gradient for CSI} and preconditioned conjugate gradient algorithms without
decreasing the solution quality. Indeed, for both algorithms, convergence is about 3 time faster when $\lambda$ is set by hand while the mean square errors are almost the same.

Finally, preconditioned conjugate gradient is from 5 to 10 times faster than \emph{alternated conjugate gradient for CSI} according to the chosen stopping rule and the value of $\lambda$. Further research should be conducted in order to explain
such large variations in convergence speed, but a first analysis suggests that the preconditioned conjugate gradient algorithm is particularly efficient for "easy" object under test, \ie for small object under test and/or object under tests
with a low contrast.

We performed the same two tests with the \emph{two cylinder} object under test using the same parameters than for the \emph{one cylinder} case. Results were quite similar and are not
presented here. They nevertheless confirm, as  stated in Subsection \ref{subsec:facteur_poid}, that a set of parameters $\lambda$ and
$\lambda\reg$ which is efficient for a given contrast, will remain efficient for a whole set of \emph{similar} contrasts. This then confirms the
possibility of choosing the value of $\lambda$ and $\lambda\reg$ by performing a training step using known contrasts.

\section{Conclusion}

Both the criterion form and the optimization scheme of the CSI method were analyzed. We established that the weight factor prescribed in the CSI
method was not suitable and could lead to a degenerate solution. We also underlined that the CSI optimization scheme does not take advantage of
conjugacy in an optimal way.

We then proposed two new methods, both making use of the same criterion which is similar to the one used in CSI. However, the weight factor is
set heuristically. We put forward that solution quality and computation time were directly related to the value of this factor. The role of the
regularization term was also investigated and we chose to use a regularization penalty term in our criterion.

\emph{Alternated conjugate gradient for CSI}, despite nested iterative algorithms, appears to be faster than CSI. We also proposed a preconditioned conjugate gradient algorithm for simultaneous
updates of the unknowns. The latter scheme is insensitive to the relative scale of the data. Our results indicate that preconditioned conjugate gradient is sometimes
faster and sometimes slower than block component approaches. Indeed, their respective behavior strongly depends  on the nature of object under test.

A more precise study of the behavior of block-component and simultaneous update approaches should be undertaken in order to better evaluate
their relative performance. We could then expect to determine which algorithm should be favored according to the experimental conditions.

Obviously, an unsupervised method for tuning the weight factor could be of great interest. In some cases the hand-tuning of $\lambda$ from
training step is a good choice but it may be sometimes difficult to apply.

\appendix
\subsection{Proof of the degeneracy of minimizers of $F\CSI=F_1+\lambda\CSI F_2$}
\label{degenerate}

Let us first assume that there exist a set \Wv of current distributions $\wb_i$ such that $F_1$ cancels, \ie
\begin{equation}
\label{eq:cte1} \yb_i=\Gv_\oD\wb_i
\end{equation}
and a point of the domain $D$ where the total field $\Eb^0_i+\Gv_\cD\wb_i$ simultaneously cancels for all $i$, \ie
\begin{equation}
\label{eq:vi} \eb_k^T(\Eb^0_i+\Gv_\cD\wb_i)=0
\end{equation}
for some $k\in\acc{1\ldotsv n}$ where $\eb_k$ is the $k$th basis vector. Then we have
\begin{align}
&F(\xb,\Wv)\to F_1 =0
\end{align}
provided that $\lambda\reg=0$ and $x_k\to\infty$. Indeed, in such conditions, $F_2$ does not depend on $x_k$, according to~\eqref{eq:F2}
and~\eqref{eq:vi}, while the denominator of $\lambda\CSI$ does, in such a way that $\lambda\CSI$ decreases to zero when $x_k$ takes arbitrarily
large values.

Finally, let us justify the existence of such a set of current distributions. For fixed values of $i$ and $k$, constraints~\eqref{eq:cte1}
and~\eqref{eq:vi} together form a linear system of $N+1$ equations, depending of the $n$ unknowns formed by $\wb_i$. Since $n\gg N$ (typically,
$n\approx N^2$), this system is likely to be undeterminate, \ie it has infinitely many solutions. This actually holds for any values of $i$ and
$k$, hence the proof.


\subsection{Conjugate gradient algorithms for complex unknown values}
\label{CG} We present here the linear and nonlinear Polak-Ribière-Polyak conjugate gradient algorithms in Table \ref{algo:LCG_general} and Table
\ref{algo:NLCG_general}, respectively . In all cases, $\gb=\nabla_\vb f$ denotes the gradient of a criterion $f$ \wrt an unknown vector $\vb$,
$\pb$ the descent direction induced by conjugate gradient, $\beta$ the conjugacy factor that allows to compute the successive conjugate descent directions $-\gb
+ \beta \pb$, $\alpha$ the optimal step length $\argmin_{\alpha}f(\xb+\alpha \pb^k)$, and $k$ the current iteration number.

Tables~\ref{algo:LCG_general} and \ref{algo:NLCG_general} address the case of complex-valued unknown vectors, which corresponds to the relevant
situation here, while classical numerical optimization textbooks such as \cite{Nocedal_livre_NO} are restricted to real-valued vectors. Indeed,
it can be shown that our complex-valued implementations of conjugate gradient produce the same iterations as the reference real-valued schemes, where all
complex-valued quantities would be replaced by an equivalent real-valued representation. For instance, \xb should be substituted by a real
vector of length $2n$ such as $[\Re(\xb)\T,\,\Im(\xb)\T]\T$.

\begin{table}[t]
\setlength\fboxsep{2mm}
\begin{spacing}{1.2}
\fbox{\kern-8pt\begin{minipage}{.985\linewidth}
\smallskip
\begin{algorithmic}
\STATE\COMMENT{Minimization of a quadratic function $f(\vb)=\vb^\dag\Av\vb-2\Re(\bb^\dag\vb)+c$ \wrt $\vb$, where \Av is a Hermitian matrix, \vb
and \bb complex-valued vectors and $c$ a real-valued constant} \STATE Initialize \vb \STATE $\gb\gets\Av\vb-\bb$ \COMMENT{\ie $\gb=\nabla_\vb
f/2$} \STATE $k \gets 0$ \REPEAT \STATE $\rho=\stdnorm{\gb}^2$ \IF{$k=0$} \STATE $\pb\gets-\gb$ \ELSE \STATE $\beta\gets\froc{\rho}{\rho\old}$
\STATE $\pb\gets-\gb+\beta\pb$ \ENDIF \STATE $\hb\gets\Av\pb$ \STATE $\alpha\gets\froc{\rho}{\pb^\dag\hb}$ \STATE $\vb\gets\vb+\alpha \pb$
\STATE $\gb\gets\gb+\alpha \hb$ \STATE $\rho\old\gets\rho$ \STATE $k \gets k+1$ \UNTIL{Sufficient decrease of $\rho=\norm{\gb}^2$}
\end{algorithmic}
\end{minipage}}
\end{spacing}
\caption{Linear conjugate gradient algorithm in the case of a complex-valued unknown vector} \label{algo:LCG_general}
\end{table}

\begin{table}[t]
\setlength\fboxsep{2mm}
\begin{spacing}{1.2}
\fbox{\kern-8pt\begin{minipage}{.985\linewidth}
\smallskip
\begin{algorithmic}
\STATE\COMMENT{Minimization of a non quadratic function $f$ \wrt $\vb$} \STATE Initialize \vb \STATE $k \gets 0$ \REPEAT \STATE
$\gb\gets\nabla_\vb f$ \IF{$k=0$} \STATE $\pb\gets-\gb$ \ELSE \STATE $\beta\gets\froc{\Re\pth{(\gb-\gb\old)^\dag\gb}}{\stdnorm{\gb\old}^2}$
\STATE $\pb\gets-\gb+\beta\pb$ \ENDIF \STATE $\alpha\gets\argmin_{\alpha}f(\vb+\alpha \pb)$ \STATE $\vb\gets\vb+\alpha \pb$ \STATE
$\gb\old\gets\gb$ \STATE $k \gets k+1$
\UNTIL{Sufficient decrease of $\norm{\gb}^2$}
\end{algorithmic}
\end{minipage}}
\end{spacing}
\caption{Nonlinear Polak-Ribière-Polyak conjugate gradient algorithm in the case of a complex-valued unknown vector} \label{algo:NLCG_general}
\end{table}

\subsection{Optimal step length calculation for simultaneous optimization schemes}
\label{third}

For the simultaneous optimization schemes, the optimal step length $\hat{\alpha}$ is the minimizer of
$F(\xb+\alpha\pb_x,[\wb_i+\alpha\pb_i]_i)$.
According to~\eqref{eq:2termes_typique},~\eqref{eq:F1},~\eqref{eq:F2} and~\eqref{eq:Freg_matrice}, $F$ is a quartic polynomial function of
$\alpha$ with real coefficients
\begin{equation}
\label{eq:simul_polynome} F(\alpha)=R_0+\alpha R_1+\alpha^2R_2+\alpha^3R_3+\alpha^4R_4
\end{equation}
where
\begin{align}
\label{eq:def_poly_maj_simul}
R_0=&F(\xb,\Wv)\notag\\
R_1=&2\Re\left(-\sum_i(\rb_{1i}^\dag\Go\pb_i+\lambda\rb_{2i}^\dag\qb_{2i})+\lambda\reg\rb\reg^\dag\pb_x\right)\notag\\
R_2=&\sum_i\left(\|\Go\pb_i\|^2+ \lambda\left(2\Re(\rb_{2i}^\dag \sb_{2i})+\|\qb_{2i}\|^2\right)\right)\notag\notag\\
&+ \lambda\reg\|\Dv\pb_x\|^2\notag\\
R_3=&2\lambda\sum_i\Re(\qb_{2i}^\dag \sb_{2i})\notag\\
R_4=&\lambda\sum_i\|\sb_{2i}\|^2\notag
\end{align}
and
\begin{align}
&\rb_{1i}=\yb_i-\Go\wb_i\notag\\
&\rb_{2i}=\Xv(\Eb^0_i+\Gc\wb_i)-\wb_i\notag\\
&\rb\reg=\Dv\xb\notag\\
&\qb_{2i}=\diag{\pb_x}(\Eb^0_i+\Gc\wb_i)+(\Xv\Gc-\Iv)\pb_i\notag\\
&\sb_{2i}=\diag{\pb_x}\Gc\pb_i\ \ . \notag
\end{align}
By necessary condition, the minimizer $\hat{\alpha}$ cancels the derivative of $F(\alpha)$
\begin{equation}
\label{eq:deriv_F_simul} F'(\alpha)=R_1+2\alpha R_2+3\alpha^2R_3+4\alpha^3R_4
\end{equation}
which is a cubic polynomial function with real coefficients. Thus, $\hat{\alpha}$ can be calculated exactly: it suffices to determine the
solutions of $F'(\alpha)=0$ (for instance, by Cardano's method), and to select the one that minimizes $F(\alpha)$ among the real solutions.

\subsection{Expression of preconditioner $\Pv$}
\label{fourth} We detail the expression of the proposed preconditioner $\Pv$. In the preconditioned conjugate gradient method the unknown vector is formed by $(\xb,\Wv)$. The
proposed preconditioner is then the inverse of the diagonal matrix formed by the diagonal entries of $\Qv$ followed by $M$ instances of the
diagonal entries of $\Av$
\begin{equation}
  \label{eq:P}
  \Pv=\diag{\diag{\Qv}^\tD,\diag{\Av}^\tD\ldotsv\diag{\Av}^\tD}^{-1}\ \ .
\end{equation}
According to \eqref{eq:def_A} and \eqref{eq:def_Q}, all these elements can be computed off-line or at a relatively low cost (term by term vector
multiplications) except for $\diag{(\Xv\Gv_\cD-\Iv)^\dag(\Xv\Gv_\cD-\Iv)}$ that can be expressed as
\begin{align}
\diag{(\Xv\Gv_\cD-\Iv&)^\dag(\Xv\Gv_\cD-\Iv)}\notag\\
&=|\Gv_\cD|^2|\xb|^2-2\Re(\diag{\diag{\Gv_\cD}}\xb)-1\notag
\end{align}
where $|\Gv_\cD|^2$ and $|\xb|^2$ represent the matrix and the vector formed with the square modulus of the entries of $\Gv_\cD$ and $\xb$, respectively. Therefore, evaluation of this term essentially requires one multiplication by a $n \times n$ full matrix at each iteration.

\section*{Acknowledgment}
Support for this work was provided under programs: National Sciences and Engineering Research Council of Canada (NSERC) Ph.D. Grant, Fonds
québécois de la recherche sur la nature et les technologies (FQRNT) Ph.D. Grant, ANR-OPTIMED, FQNRT Grant \#PR-113254 and NSERC Discovery Grant
\#138417-06.
\bibliographystyle{IEEEtran}
\bibliography{bibenabr,revuedef,revueabr,TMO_ref_gen}

\end{document}